\documentclass[a4paper, 11pt]{amsart}

\usepackage[latin1]{inputenc}
\usepackage{amssymb, amsmath}
\usepackage{mathrsfs}
\usepackage{amsthm}
\usepackage[english]{babel}

\usepackage[numeric, lite, initials, nobysame]{amsrefs}
\input xy
\xyoption{arrow} \xyoption{matrix} \xyoption{frame}
\xyoption{curve}

\let\oldmarginpar\marginpar
\renewcommand\marginpar[1]{\-\oldmarginpar[\raggedleft\footnotesize #1]%
{\raggedright\footnotesize #1}}

\def\:{\colon}
\def\<{\left\langle}
\def\<{\right\rangle}
\def\({\left(}
\def\){\right)}
\def\epsilon{\varepsilon}
\def\phi{\varphi}
\def\subset{\subseteq}

\def\leq{\leqslant}
\def\geq{\geqslant}

\def\lra{\longrightarrow}

\def\mapsto{\longmapsto}
\def\smash{\wedge}
\def\xra{\xrightarrow}

\def\op{\text{op}}

\def\wh{\widehat}
\def\odd{\text{odd}}

\def\DR{\mathscr{D}_R}
\def\Cl{\mathcal C \ell}
\def\SS{\mathbb S}\def\C{\mathbb C}

\def\Z{\mathbb Z}
\def\F{\mathbb F}

\newtheorem{thm}{Theorem}
\newtheorem{lem}[thm]{Lemma}
\newtheorem{prop}[thm]{Proposition}
\newtheorem{cor}[thm]{Corollary}

\newtheorem{thmi}{Theorem}

\newtheorem{propi}[thmi]{Proposition}
\newtheorem{cori}[thmi]{Corollary}

\newtheorem*{thm*}{Theorem}
\newtheorem*{cor*}{Corollary}
\newtheorem*{prop*}{Proposition}

\numberwithin{equation}{section}
\numberwithin{thm}{section}

\theoremstyle{remark}
\newtheorem{rem}[thm]{Remark}

\theoremstyle{definition}
\newtheorem{defn}[thm]{Definition}

\newtheorem{nota}[thm]{Notation}

\DeclareMathOperator{\Hom}{Hom}\DeclareMathOperator{\End}{End}

\DeclareMathOperator{\hocolim}{hocolim}

\DeclareMathOperator{\hatotimes}{\widehat\otimes}
\DeclareMathOperator{\DDer}{\mathscr Der}

\DeclareMathOperator{\Der}{Der}\DeclareMathOperator{\Prod}{Prod}
\DeclareMathOperator{\Bil}{Bil}
 \DeclareMathOperator{\Sym}{Sym}
\DeclareMathOperator{\Alt}{Alt} \DeclareMathOperator{\Quad}{QF}
\DeclareMathOperator{\Asym}{Asym} \DeclareMathOperator{\per}{Per}

\newcommand{\ie}{i.e.}

\newcommand{\eg}{e.g.}

\title{Quadratic forms classify products on quotient ring spectra }

\author{A. Jeanneret}
\address{Mathematisches Institut, Sidlerstrasse 5, 3012 BERN, Switzerland}
\email{alain.jeanneret@math.unibe.ch}

\author{S. W\"uthrich}
\address{SBB, Br\"uckfeldstrasse 16, 3000 BERN, Switzerland}
\email{samuel.wuethrich@sbb.ch}


\date{15.\@ March 2010}

\begin{document}

\begin{abstract}
We construct a free and transitive action of the group of bilinear
forms $\Bil(I/I^2[1])$ on the set of $R$-products on $F$, a
regular quotient of an $E_\infty$-ring spectrum $R$ with $F_*
\cong R_*/I$. We show that this action induces a free and
transitive action of the group of quadratic forms
$\Quad(I/I^2[1])$ on the set of equivalence classes of
$R$-products on $F$. The characteristic bilinear form of $F$
introduced by the authors in a previous paper is the natural
obstruction to commutativity of $F$. We discuss the examples of
the Morava $K$-theories $K(n)$ and the 2-periodic Morava
$K$-theories $K_n$.

\end{abstract}

\maketitle

\section{Introduction}\label{intro}

With the advent of sound foundations for a theory of modules over
an $E_\infty$-ring spectrum $R$ (for instance as developed in
\cite{ekmm}), it has become possible to mimic in homotopy theory
well-known constructions usually performed in algebra. The setting
is the homotopy category $\DR$ of $R$-module spectra over $R$, a
category equipped with a smash product $\smash_R$ (the equivalent
of the tensor product), giving $\DR$ the structure of a symmetric
monoidal category. Objects in $\DR$ may be regarded as ordinary
spectra by neglect of structure, via a monoidal functor to the
classical stable homotopy category.

With this framework at hand, the problem of constructing quotient
spectra, \ie\ spectra whose homotopy groups are isomorphic to a
given quotient of the coefficient ring $R_*=\pi_*(R)$ of $R$,
admits a clean and transparent solution for a large class of
quotients. The quotients in question are the quotients $R_*/I$ by
ideals $I$ which are generated by regular sequences. The
$R$-module spectra realizing such quotients are often referred to
as regular quotients.

Shortly after the publication of \cite{ekmm}, Strickland proved
that for $E_\infty$-ring spectra $R$ for which $R_*$ forms a
domain and which is trivial in odd degrees, any regular quotient
can be realized as an $R$-ring spectrum, \ie\ as a monoid in
$\DR$, and therefore in particular as a ring spectrum
\cite{strickland}.

The aim of the present article is to give a conceptual description of the set of all
$R$-ring structures on regular quotients $F$ of $R$, as well as of the set of equivalence
classes of $R$-ring structures. Our result in both cases is based on a free and
transitive action of a certain abelian group canonically associated to $F$ on the set of
products.

As an application, we show that the characteristic bilinear form
$b_F$ of a regular quotient $F$, introduced by the authors in
\cite{jw}, is always symmetric and  provides a measure for the
non-commutativity of $F$.

As another application, we give a necessary and sufficient
criterion for a map of regular quotient rings $\pi\: F\to G$ to be
multiplicative, in terms of the characteristic bilinear forms.

We use our results to classify products on the $2$-periodic Morava-$K$-theories $K_n$,
which from an algebro-geometric point of view are the more natural objects to study than
their classical variants $K(n)$. In contrast to $K(n)$, we show that $K_n$ supports a
large number of products, even many commutative ones for $p$ odd and $n>1$.

In addition, we confirm many well-known facts concerning certain
families of quotients of complex cobordism $MU$, whose existing
proofs are in many cases technically forbidding and scattered in
the literature.

\smallskip
\smallskip

We now proceed to a more detailed overview of the content of this
article. Throughout, $R$ denotes an $E_\infty$-ring spectrum for
which $R_*$ is a domain and is trivial in odd degrees.

The following result assembles our two main theorems (Theorems
\ref{action} and \ref{altiso}). The symbol $I/I^2[1]$ stands for
the graded module $I/I^2$ shifted by one, where $I\subseteq R_*$
is an ideal.

\begin{thmi}\label{maini}
Let $F$ be a regular quotient of $R$ with coefficients $F_*\cong R_*/I$.

\begin{enumerate}\itemsep2pt

\item There is a natural free and transitive action of the abelian
group $\Bil(I/I^2[1])$ of bilinear forms on $I/I^2[1]$ on the set
of $R$-products on $F$.

\item This action induces a free and transitive actions of the
abelian group $\Quad(I/I^2[1])$ of quadratic forms on $I/I^2[1]$
on the set of equivalence classes of $R$-products on $F$.

\end{enumerate}
\end{thmi}

For a regular quotient ring $F=R/I$ with product $\mu$ and a
bilinear form $\beta\in\Bil(I/I^2[1])$, we will denote by $\beta
F$ the $R$-module $F$, endowed with the product $\beta\mu$ in the
sequel.

For the proof of the theorem we build on our previous paper
\cite{jw}. The central ingredient is the module of (homotopy)
derivations $\DDer_R^*(F)$. Of crucial importance is the fact
proved in \cite{jw} that $\DDer_R^*(F)$ does not depend on the
product of $F$, as a submodule of the algebra of endomorphisms
$F^*_R(F)$.

Applied to $R=\wh{E}(n)$, the completed Johnson--Wilson theories,
and $F=K(n)$, the theorem implies immediately that for $p$ odd,
there is precisely one $\widehat{E}(n)$-product on $K(n)$, which
therefore must be commutative. For $p=2$, it implies that there
are precisely two non-equivalent $\widehat{E}(n)$-products on
$K(n)$. They are both non-commutative, as we will see below. These
are well-known results.

For $R=E_n$, the Morava $E$-theories, and $F=K_n$, we deduce that
there are $p^n n^2$ different $E_n$-products and $p^n \frac{n}{2}
(n+1)$ equivalence classes of $E_n$-products on $K_n$. By
construction, it follows that all the products remain different
when regarded as products on the underlying spectra $K_n$.

It is natural to ask whether there is an invariant which distinguishes the different
products on $F$ or at least the different equivalence classes of products. A candidate is
the characteristic bilinear form
\[
b_F\: I/I^2[1] \otimes_{F_*} I/I^2[1] \lra F_*
\]
of a regular quotient ring $F$ constructed in \cite{jw}. We prove as Corollary
\ref{isoquad}:

\begin{propi}
The characteristic bilinear forms of equivalent products on $F$
coincide. The converse holds whenever $F_*$ is $2$-torsion free.
\end{propi}

In fact, the characteristic bilinear form $b_F$ admits a natural
characterization in terms of the action of the theorem. To express
it, let $F^\op$ denote the opposite ring of $F$. We prove as
Corollary \ref{relateffop}:

\begin{propi}\label{charbfi}
The characteristic bilinear form $b_F$ of a regular quotient ring $F$ satisfies $F^\op =
b_F F$.
\end{propi}

Hence $b_F$ is the obstruction to commutativity of $F$:

\begin{cori}\label{commcriti}
A regular quotient ring $F$ is commutative if and only if $b_F=0$.
\end{cori}

Consider again the Morava $K$-theories $K(n)$ at $p=2$. We proved
in \cite{jw} that it admits an $\wh{E}(n)$-product $\mu$ with
non-trivial characteristic bilinear form. Corollary
\ref{commcriti} implies that $\mu$ cannot be commutative.
Therefore the second product on $K(n)$ is neither, as it must be
the opposite of $\mu$. Moreover, Proposition \ref{charbfi}
recovers the well-known formula (see Section \ref{examples} for
the definition of $v_n$ and $Q_{n-1}$)
\[
\mu^\op= \mu \circ (1+ v_n Q_{n-1} \wedge Q_{n-1}).
\]

As a consequence of Theorem \ref{maini}, there is in general a
large variety of products on $F$, even up to equivalence, unless
there are only few bilinear forms on $I/I^2[1]$ due to sparseness
of the coefficients $F_*$. One may ask if the situation changes
when one restricts to commutative products. To approach this
question, one needs a formula which expresses how $b_F$ transforms
under the action of $\Bil(I/I^2[1])$, in view of Corollary
\ref{commcriti}. Let $\beta^t$ denote the transpose of a bilinear
form $\beta$ on $I/I^2[1]$, defined by $\beta^t(x\otimes y) =
\beta (y\otimes x)$. We prove as Corollary \ref{newcharbil}:

\begin{propi}\label{commclassi}
Let $F$ be a regular quotient ring with characteristic bilinear
form $b_F$ and let $\beta$ be a  bilinear form in
$\Bil(I/I^2[1])$. Then the characteristic bilinear form of $\beta
F$ is given by $b_{\beta F} = b_F - \beta -\beta^t$.
\end{propi}

With Corollary \ref{commcriti}, it follows that for commutative
$F$, $\beta F$ is commutative if and only if $\beta$ is
antisymmetric. Together with Theorem \ref{maini}, this implies the
following result (Proposition \ref{comprod}, Corollary
\ref{commprodequiv}), which sharpens a result of
\cite{strickland}.

\begin{cori}
Let $F$ be a regular quotient ring of $R$. If $2\in F_*$ is invertible, there exists a
unique commutative product on $F$ up to equivalence. If $F_*$ is $2$-torsion free, there
exists at most one commutative product on $F$ up to equivalence.
\end{cori}

For the $2$-periodic Morava $K$-theories, Proposition
\ref{commclassi} implies that there are $p^n\frac{n}{2}(n-1)$
commutative $E_n$-products for odd $p$, all of which are
equivalent. At the prime $2$, $K_n$ admits a product with
non-trivial characteristic bilinear form, by a result from
\cite{jw}. From this, it follows that there does not exists any
commutative product on $K_n$ for $p=2$.

Using the fact proved in \cite{strickland} that the Brown--Peterson spectrum $BP$ at a
prime $p$ admits a commutative $MU$-product, it follows that there is a unique
commutative $MU$-product on $BP$ up to equivalence.

The products on regular quotients $F$ constructed in
\cite{strickland} have a very special form. To explain in what
sense, let $(x_1, x_2, \ldots)$ be a regular sequence generating
$I$, where $F_*\cong R_*/I$. Then $F$ is equivalent as an
$R$-module spectrum to $R/x_1\smash_R R/x_2\smash_R \cdots$ (see
Section \ref{recollection} for details). The products considered
in \cite{strickland} are all obtained by ``smashing together''
products on the $R$-module spectra $R/x_k$. We call such products
diagonal and refer to products equivalent to diagonal ones as
diagonalizable. In \cite{jw}, we showed that the characteristic
bilinear form of a diagonal regular quotient ring is diagonal.
Together with Proposition \ref{commclassi}, this implies
(Corollary \ref{symetric}):

\begin{cori}
The characteristic bilinear form $b_F$ of a regular quotient ring $F$ is symmetric.
\end{cori}

The following result is proved as Proposition \ref{diagprod}:

\begin{propi}
Assume that $R_*$ is a finite-dimensional regular local ring with maximal ideal $I$ and
suppose that $F$ is an $R$-ring satisfying $F_*\cong R_*/I$. If the characteristic $p$ of
$F_*$ is zero or an odd prime, then $F$ is diagonalizable. If $p=2$, then $F$ is
diagonalizable unless $b_F$ is alternating and non-trivial, in which case $F$ is not
diagonalizable.
\end{propi}

This implies for instance that any $E_n$-product on $K_n$ is diagonalizable, for $p$
arbitrary. However, not every regular quotient ring is diagonalizable: We construct a
non-diagonalizable $MU$-ring spectrum in Section \ref{examples}.

As an application of Theorem \ref{maini}, we give a necessary and sufficient condition
for a map $\pi\: F\to G$ between regular quotients of $R$ to be multiplicative. Let
$I\subset J$ be the ideals of $R_*$ for which $F_*\cong R_*/I$ and $G_*=R_*/J$,
respectively. In \cite{jw}, we introduced a bilinear form
\[
b^G_F\: (G_*\otimes_{F_*} I/I^2[1]) \otimes_{G_*} (G_*\otimes_{F_*} I/I^2[1]) \to G_*,
\]
which depends on $\pi$. Let $b_F$ and $b_G$ denote the
characteristic bilinear forms of $F$ and $G$, respectively. Let
$\pi^*(b_G)$ be the bilinear form on $G_*\otimes_{F_*} I/I^2[1]$
obtained by ``pulling back'' $b_G$ along the morphism $\bar\pi \:
I/I^2[1] \to J/J^2[1]$ induced by $\pi$.

\begin{thmi}
Suppose that $\pi\: F\to G$ is as above and assume that the
induced map $G_*\otimes_{F_*} I/I^2[1] \to J/J^2[1]$ is injective.
Then  $\pi$ is multiplicative if and only if $G_* \otimes
b_{F}=b^{G}_{F}=\pi^*(b_G)$.
\end{thmi}

As an illustration, we show that there are infinitely many $MU$-products on the spectrum
$P(n)$ for any prime $p$ such that the canonical map $BP \to P(n)$ is multiplicative,
where $BP$ is endowed with an arbitrary commutative $MU$-product (see Section
\ref{examples}).

\subsection*{Relation to other work}
The proof of Theorem \ref{maini} requires a formula stated in
\cite{ang}, which gives a description of the set of all products
on a regular quotient of $R$. This formula may also be viewed as
providing an answer to the question of how to classify products on
regular quotients. However, its technical formulation basically
forbids any serious practical application. Because the proof given
in \cite{ang} appears rather incomplete and fragmentary, we give
an independent and complete proof here.

\subsection*{Acknowledgments}
The second author would like to thank Prof.\@ Kathryn Hess for her support throughout his
time at the EPFL in Lausanne.

\subsection*{Notation and conventions}

In this article, we will work in the framework of $\SS$-modules of \cite{ekmm}. In this
setting, $E_\infty$-ring spectra correspond to commutative $\SS$-algebras. Throughout,
$R$ denotes an even commutative $\SS$-algebra, \ie\ one with $R_{\odd} =0$. We also
assume that the coefficient ring $R_*$ of $R$ is a domain (see \cite{jw}*{Remark 2.11}).
Associated to $R$ is the homotopy category $\DR$ of $R$-module spectra. For simplicity,
we refer to its objects as $R$-modules. The smash product $\smash_R$ endows $\DR$ with a
symmetric monoidal structure. We will abbreviate $\smash_R$ by $\smash$ throughout the
paper.

Monoids in $\DR$ are called $R$-ring spectra or just $R$-rings. Unless otherwise
specified, we use the generic notation $\eta_F\: R\to F$ (or simply $\eta$) for the unit
and $\mu_F\: F\smash F\to F$ (or simply $\mu$) for the multiplication of an $R$-ring $F$.
Mostly, $\eta_F$ will be clear from the context, in which case we call a map $\mu_F\:
F\smash F\to F$ which gives $F$ the structure of an $R$-ring an $R$-product or just a
product. For a given $R$-ring $(F, \mu_F, \eta _F)$, we will often be in the situation
where we consider another product $\bar \mu_F$ on $F$. We then write $\bar F$ for the
$R$-ring $(F, \bar \mu_F, \eta _F)$. We denote the opposite of an $R$-ring $F$ by
$F^\op$. Its product is given by $\mu_{F^\op}= \mu_F\circ\tau$, where $\tau \: F\smash
F\to F\smash F$ is the switch map.

An $R$-ring $(F, \mu_F,\eta_F)$ determines multiplicative homology
and cohomology theories $F_*^R(-)=\pi_*(F\smash -) =
\DR^{-*}(R,F\smash -)$ and $F^*_R(-)=\DR^*(-, F)$, respectively,
on $\DR$. For an $R$-module $M$, the homology $F_*^R(M)$ is an
$F_*$-bimodule in a natural way. Even if $F_*$ is commutative, the
left and right $F_*$-actions may well be different. However, if we
assume that $F$ is a quotient of $R$, by which we mean that the
unit map $\eta_F$ induces a surjection on homotopy groups (see
Section 2 below for definitions), the left and right $F_*$-actions
agree. In this case, we can refer to $F_*^R(M)$ as a $F_*$-module
without any ambiguity. A similar discussion applies to cohomology
$F_R^*(M)$. See Section 1.1 of \cite{jw} for a more detailed
discussion.

We write $M_*[d]$ for the $d$-fold suspension of a graded abelian
group $M_*$, so $(M_*[d])_k = M_{k-d}$. With this convention, we
have $(\Sigma^d M)_* = M_*[d]$ for an $R$-module $M$. We use the
convention $M^*=M_{-*}$. If the ground ring is clear from the
context, we omit it from the tensor product symbol $\otimes$ from
now on. We write $D_{F_*}(M_*)$ or just $D(M_*)$ for the dual
$\Hom^*_{F_*}(M_*, F_*)$ of a graded module $M_*$ over a graded
ring $F_*$.

We introduce some notation and recall some well-known facts
concerning bilinear and quadratic forms. For an $F_*$-module $V$,
we write $\Bil(V)$ for the abelian group of (degree 0) bilinear
forms on $V$. For $\beta \in \Bil(V)$, we set $\beta^t(x\otimes
y)= \beta(y \otimes x)$ for $x,y \in V$. A bilinear form $\beta
\in \Bil(V)$ is {\it symmetric} if $\beta^t= \beta$, {\it
antisymmetric} if $\beta^t= -\beta$ and {\it alternating} if
$\beta(v\otimes v)=0$ for any $v \in V$. We write $\Sym(V)$,
$\Asym(V)$ and $\Alt(V)$ for the subgroups of $\Bil(V)$ consisting
of the symmetric, antisymmetric and alternating bilinear forms,
respectively. If $V$ is 2-torsion-free, we have $\Sym(V) \cap
\Alt(V) = 0$ and $\Asym(V)=\Alt(V)$. If $2 \in F_*$ is invertible,
we have the usual decomposition $\Bil(V)=\Sym(V) \oplus \Alt(V)$.

Let $\Quad(V)$ denote the group of quadratic forms $q \: V \to
F_*$. Recall that the grading convention is that $|q(v)|=2n$ for
$v\in V_n$. For $\beta \in \Bil(V)$, $q(v)= \beta (v \otimes v)$
is easily seen to be a quadratic form. We thus obtain a group
homomorphism $\chi \: \Bil(V)\to \Quad(V)$, whose kernel is
$\Alt(V)$. If $V$ is $F_*$-free, $\chi$ is surjective (see
\cite{bourbaki}*{Chap.\@ IX, §3, Prop.\@ 2}) and so we have a
canonical isomorphism $\Bil(V)/\Alt(V) \cong \Quad(V)$. If $2 \in
F_*$ is invertible, we recover the well-known isomorphism $\Sym(V)
\cong \Quad(V)$.

For a ring homomorphism $\pi_* \: F_* \to k_*$ and $\beta \in
\Bil(V)$, we define $k_*\otimes \beta$ to be the bilinear form on
the $k_*$-module $k_*\otimes_{F_*} V$ determined by
\[
(k_*\otimes \beta)((1\otimes x) \otimes (1\otimes y)) = \pi_*(\beta( x\otimes y)),
\]
for $x,y\in V$. If $\pi_* \: W \to V$ is a morphism of
$F_*$-modules and $\beta \in \Bil(V)$, $\pi^*(\beta)$ denotes the
bilinear form on $W$ which on $x,y\in W$ takes the value
\[
\pi^*(\beta)(x \otimes y) = \beta(\pi_*( x)\otimes \pi_*(y)).
\]

\section{Recollection}\label{recollection}
In this section we collect some results, definitions and notation from \cite{jw} which we
are using in the present paper.

A {\em quotient module}\/ of $R$ is an $R$-module $F$ with a map
of $R$-modules $\eta_F\: R\to F$ which induces a surjection on
homotopy groups, that is $F_* \cong R_*/I$. We will write $F=R/I$
for such an $F$ in the sequel. The modules of interest for our
purposes are the {\em regular quotient modules}\/ of $R$. By this,
we mean quotient modules $F=R/I$ whose ideal $I$ is generated by
some (finite or infinite) regular sequence $(x_1, x_2, \ldots)$ in
$R_*$.

A (regular) {\em quotient ring}\/ of $R$ is an $R$-ring $(F,
\mu_F, \eta_F)$ with product $\mu_F$ such that $(F, \eta_F)$ is a
(regular) quotient module of $R$. For instance, let $F=R/I$ be a
regular quotient of $R$ and $(x_1,x_2,\ldots)$ a regular sequence
generating the ideal $I$. Then $F$ is isomorphic in $\DR$ to
\[
R/x_1\smash R/x_2 \smash \cdots := \hocolim_k R/x_1\smash \cdots
\smash R/x_k,
\]
where for $x\in R_*$, we denote by $R/x$ the homotopy cofibre of
$x\: \Sigma^{|x|}R \to R$. For any products $\mu_i$ on $R/x_i$,
there is a uniquely determined product $\mu$ on $F=R/I$ such that
the natural maps $j_i \: R/x_i \to F$ are multiplicative and
commute for $k\neq l$, \ie\ $\mu (j_k\smash j_l) = \mu^\op
(j_k\smash j_l)$. This ring $F$ is called the {\em smash ring
spectrum}\/ of the $R/x_i$. If we need to be more precise, we
refer to the product map $\mu_F$ as the {\it smash ring product}\/
of the $\mu_i$. A regular quotient ring $F$ whose product is of
this form is said to be {\em diagonal} or {\em diagonal with
respect to $(x_1,x_2,\ldots)$} if we need to keep track of the
regular sequence.

An {\it admissible pair} is a triple $(F, k, \pi)$ consisting of two quotient $R$-rings
$(F, \mu_F, \eta_F)$, $(k, \mu_k, \eta_k)$ and a unital $R$-module map $\pi\: F\to k$,
\ie\ an $R$-morphism $\pi$ with $\pi \eta_F = \eta_k$. If $\pi$ is a map of $R$-ring
spectra, we call $(F, k, \pi)$ a {\em multiplicative admissible pair}. A typical example
of an admissible pair is $(F,F,1_F)$ where $1_F$ is the identity on $F$, but where we
distinguish two products $\mu$ and $\nu$ on $F$.

In the following, we fix an admissible pair $(F=R/I, k, \pi)$. Its
{\em characteristic homomorphism} is a homomorphism of
$F_*$-modules
\begin{equation}\label{natmap}
\phi^k_F\:  I/I^2[1] \lra k_*^R(F),
\end{equation}
which is natural in $F$ and $k$ and independent of the products on $F$ and $k$.

The homology group $k_*^R(F)$ carries a natural $k_*$-algebra
structure, whose product is defined by the following composition
of $k_*$-homomorphisms
\begin{equation}\label{multalgebra}
m^{k}_{F} \: k_*^R(F) \otimes_{k_*} k_*^R(F) \xra{\kappa_{k}} k_*^R(F\smash F).
\xra{k_*^R(\mu_F)} k_*^R(F),
\end{equation}
Here $\kappa_{k}$ stands for the K\"unneth homomorphism associated to the ring $k$.

The {\em characteristic bilinear form} $b^k_F$ associated to
$(F,k,\pi)$ is defined as the following composition of
$k_*$-homomorphisms
\[
b^k_F\: (k_* \otimes_{F_*} I/I^2[1])^{\otimes 2} \xra{{}^\pi\phi^{\otimes 2}}
k^R_*(F)^{\otimes 2} \xra{m^k_F} k_*^R(F) \xra{k_*(\pi)}k_*^R(k) \xra{(\mu_k)_*}k_*,
\]
where ${}^\pi\phi$ is the $k_*$-homomorphism canonically induced by $\phi$. The {\em
characteristic quadratic form} $q^k_F$ is defined as $q^k_F( \bar x)=b^k_F(\bar x\otimes
\bar x)$ for $\bar x \in k_* \otimes_{F_*} I/I^2[1]$.

We write $\phi_F$, $b_F$ and $q_F$ for the characteristic
homomorphism, bilinear and quadratic forms of the admissible pair
$(F,F, 1_F)$, respectively, for a quotient ring $F$.

If $(F,k,\pi)$ is multiplicative, the characteristic bilinear form
$b^k_{F^\op}$ of the associated admissible pair $(F^\op, k, \pi)$
is trivial. In particular, $b^F_{F^\op}=0$ for a quotient ring $F$
and $b_F=0$ for a commutative quotient ring $F$.

The characteristic homomorphism $\phi=\phi^k_F$ lifts to an
algebra homomorphism
\[
\Phi\:  \Cl(k_*\otimes_{F_*} I/I^2[1], q^k_F) \lra k^R_*(F),
\]
where $\Cl(k_*\otimes_{F_*} I/I^2[1], q^k_F)$ denotes the Clifford algebra of the
quadratic module $(k_*\otimes_{F_*} I/I^2[1], q^k_F)$. If $F$ is a regular quotient, then
$\Phi$ is an isomorphism. In particular, this yields an algebra isomorphism
\[
F_*^R(F^\op) \cong \Lambda(I/I^2[1]).
\]
To be more explicit, fix a regular sequence $(x_1, x_2, \ldots)$
generating $I$. This choice determines an isomorphism
$I/I^2[1]\cong\bigoplus_{i} F_* \bar x_i$, where $\bar x_i$
denotes the residue class of $x_i$ in $I/I^2[1]$. Letting
$a_i=\phi^F_{F^\op}(\bar x_i)\in F_*^R(F)$, we have
\begin{equation}\label{ffopexplicitintro}
F_*^R(F^\op) \cong \Lambda(a_1, a_2, \ldots).
\end{equation}

If $(F,k, \pi)$ is multiplicative, we may consider the module of
(homotopy) derivations $\DDer_R^*(F, k) \subset k_R^*(F)$. By
definition, these are maps $d\: F\to \Sigma^i k$ which satisfy
$d\mu_F = \mu_k(1\smash d + d\smash 1)$. If $F=k$ and $\pi=1_F$,
we write $\DDer_R^*(F)$ instead of $\DDer_R^*(F,F)$. There is a
natural $k_*$-homomorphism
\begin{equation}\label{defpsi}
\psi\: \DDer_R^*(F, k) \to \Hom^*_{F_*}(I/I^2[1], k_*),
\end{equation}
defined by $\psi(d)(\bar x)= (\mu_k)_*k^R_*(d)(\phi^k_F(\bar x))$
for $d \in \DDer_R^*(F, k)$ and $\bar x \in I/I^2[1]$. It is a
homeomorphism if both $F$ and $k$ are regular quotient rings,
where $\DDer_R^*(F,k)$ is endowed with the subspace topology
induced by the profinite topology on $k^*_R(F)$ and
$\Hom^*_{F_*}(I/I^2[1], k_*)\cong
D_{k_*}(k_*\otimes_{F_*}I/I^2[1])$ with the dual-finite topology.
The composition
\[
\Hom^*_{F_*}(I/I^2[1], k_*) \xra{\psi^{-1}} \DDer_R^*(F,k)
\subseteq k^*_R(F)
\]
is independent of the products on $F$ and $k$. This result allows
us to construct derivations. We restrict to the case where $k=F$
is a regular quotient ring and $\pi = 1_F$ here. Let
$(x_1,x_2,\ldots)$ be a regular sequence generating the ideal $I$
and let $\bar x_i^\vee \in D_{F_*}(I/I^2[1])$ be dual to $\bar
x_i$. The Bockstein operation $Q_i \in \DDer^*_R(F)$ associated to
$x_i$ is defined by $Q_i = \psi^{-1}(\bar x_i^\vee)$.

For a regular quotient ring $F$, the inclusion $\DDer_R^*(F) \to F^*_R(F)$ lifts to a
homeomorphism of $F^*$-algebras
\begin{equation}\label{cohomologyiso}
\widehat\Lambda(\DDer_R^*(F))\cong F^*_R(F),
\end{equation}
where $\widehat\Lambda(\DDer_R^*(F))$ denotes the completed exterior algebra on
$\DDer_R^*(F)$ and where $F^*_R(F)$ is endowed with the profinite topology.

\section{The action of bilinear forms on products}\label{def of action}

In this section, we show that there is a canonical action of the
group of bilinear forms  $\Bil(I/I^2[1])$ on the set of products
on a regular quotient $F=R/I$.

Let $F=R/I$ be a regular quotient and let $\Prod_R(F)\subset
F_R^*(F\smash F)$ denote the set of all products on $F$. Let
$\per(\Prod_R(F))$ be the group of permutations of the set
$\Prod_R(F)$.

Writing $V$ for $I/I^2[1]$, we have a linear isomorphism \cite{boardman}*{Lemma 6.15}
\begin{equation*}
\Bil(V) = D^0(V\otimes V)  \cong (D(V)\hatotimes D(V))^0.
\end{equation*}
Composing it with the isomorphism
\begin{equation*}
(D(V)\hatotimes D(V))^0 \xra{\psi^{-1}\hatotimes\psi^{-1}}
(\DDer^*_R(F)\hatotimes \DDer^*_R(F))^0
\end{equation*}
induced by the homeomorphism $\psi$ \eqref{defpsi} yields an isomorphism of $F^*$-modules
\begin{equation}\label{bilder}
\Bil(V)\cong (\DDer^*_R(F)\hatotimes \DDer^*_R(F))^0.
\end{equation}

The aim of this section is to prove the following result.

\begin{prop}\label{derconstruction}
Let $\mu\in\Prod_R(F)$ be a product and let $d, d'\in\DDer^*_R(F)$ be derivations with
$|d|=-|d'|$. Then the composition
\[
\mu_{d, d'}\: F\smash F\xra{1+d\smash d'} F\smash F\xra{\mu} F
\]
defines a product. This construction induces a group homomorphism
\begin{equation}\label{thirdstep}
\hat\pi\: (\DDer^*_R(F)\wh\otimes\DDer^*_R(F))^0 \lra
\per(\Prod_R(F)),
\end{equation}
which gives rise via \eqref{bilder} to an action of
$\Bil(I/I^2[1])$ on $\Prod_R(F)$.
\end{prop}


\begin{nota}
We refer to the action of Proposition \ref{derconstruction} as the
{\em canonical action} of $\Bil(I/I^2[1])$ on $\Prod_R(F)$ in the
sequel. The image of $(\beta, \mu) \in \Bil(I/I^2[1]) \times
\Prod_R(F)$ under the canonical action  will be denoted by $\beta
\mu$. Accordingly, $\beta F$ stands for $F$, endowed with the
product $\beta \mu$.
\end{nota}

\begin{rem}\label{coincidestrickland}
The proof of Proposition \ref{derconstruction} given below shows
that in the special case $F=R/x$, the action of
$\Bil((x)/(x)^2[1])\cong F_{2|x|+2}$ coincides with the action
defined in \cite{strickland}*{Prop.\@  3.1}.
\end{rem}

\begin{proof}[Proof of Proposition \ref{derconstruction}]
We first prove that $\bar\mu=\mu_{d,d'}$ is associative, \ie\ that
$\bar\mu ( 1 \smash \bar\mu) = \bar\mu (\bar\mu \smash 1)$. As a
consequence of the isomorphism \eqref{cohomologyiso}, derivations
anticommute, \ie\ for $d, d' \in \DDer^*_R(F)$, we have $ d d'=-d'
d$. Moreover, as a derviation, $d$ satisfies $d \mu = \mu (d
\smash 1 + 1 \smash d)$. This yields:
\[
\begin{split}
\bar\mu & (\bar\mu \smash 1)=\mu ((\mu + \mu(d \smash d')) \smash 1)+\mu(d \smash
d')((\mu + \mu(d \smash d') \smash 1)
\\
& =\mu(\mu \smash 1)( 1 + d \smash d' \smash 1 ) +\mu [ (\mu( d
\smash 1 + 1 \smash d)) \smash d' +
\\
& \qquad (\mu( d \smash 1 + 1 \smash d))( d  \smash d')  \smash d'  ]
\\
& =\mu(\mu \smash 1)[1 + d \smash d' \smash 1 + d \smash 1 \smash d' + 1 \smash  d \smash
d' - d \smash   dd' \smash d'] .
\end{split}
\]
On the other hand, we obtain:
\[
\begin{split}
\bar\mu & (1 \smash \bar\mu )=\mu (1 \smash  (\mu + \mu(d \smash d')) )+\mu(d \smash
d')(1 \smash (\mu + \mu(d \smash d')))
\\
& =  \mu(1 \smash \mu )( 1 + 1 \smash d \smash d'  ) +\mu [d \smash (\mu( d' \smash 1 + 1
\smash d'))  +
\\
& \qquad d \smash (\mu( d' \smash 1 + 1 \smash d'))( d  \smash d') ]
\\
& = \mu(1 \smash \mu )[1 + 1 \smash  d \smash  d' + d \smash d' \smash 1 + d \smash 1
\smash d' - d \smash   dd' \smash d'],
\end{split}
\]
which proves that $ \bar\mu$ is associative.

That $\bar\mu$ has $\eta_F\: R\to F$ as a two-sided unit is an easy consequence of the
fact that the composition $d \eta_F$ is trivial for a derivation $d$.

To prove that $\mu\mapsto \mu_{d,d'}$ defines a permutation of
$\Prod_R(F)$, it suffices to note that $\mu'\mapsto\mu'_{-d, d'}$
is a two-sided inverse. This follows from
\begin{equation}\label{inverseaction}
(\mu_{d,d'})_{-d, d'}= \bigl(\mu(1+d\smash d')\bigr)(1 - d\smash d') = \mu.
\end{equation}

We have shown so far that $(d, d')\mapsto \mu_{d,d'}$ defines a function
\[
\pi\: (\DDer^*_R(F)\times\DDer^*_R(F))^0 \to \per(\Prod_R(F)).
\]

Now \eqref{inverseaction} implies that $(\mu_{d, d'})_{e, e'} =
(\mu_{e, e'})_{d, d'}$ for derivations $d, d', e, e'$ with
$|d|=-|d'|$ and $|e|=-|e'|$. As a consequence, $\pi$ factors as
\begin{equation}\label{ccomm}
(\DDer^*_R(F)\times\DDer^*_R(F))^0 \xra{\pi'} C \subseteq
\per(\Prod_R(F)),
\end{equation}
where $C$ denotes the centre of $\per(\Prod_R(F))$. Using the
facts that i) \nolinebreak deri\-vations square to zero and ii)
that $F$ is a quotient ring of $R$, one checks that $\pi'$ in
\eqref{ccomm} is bilinear. Hence $\pi$ induces a group
homomorphism
\[
\bar\pi\: (\DDer^*_R(F)\otimes\DDer^*_R(F))^0 \lra
\per(\Prod_R(F)).
\]

Recall that $\DDer_R^*(F)$ carries the topology inherited by the profinite topo\-logy on
$F^*_R(F)$. We now show that $\bar\pi$ lifts to a group homomorphism
\[
\hat\pi\: (\DDer^*_R(F)\wh\otimes\DDer^*_R(F))^0 \lra
\per(\Prod_R(F)).
\]

Let $\End_R^*(F\smash F)$ denote $(F\smash F)^*_R(F\smash F)$.
Consider the homomorphism of monoids --- with respect to addition
and composition, respectively---
\[
\alpha\: (\DDer^*_R(F)\otimes\DDer^*_R(F))^0 \to \End_R^*(F\smash F),
\]
given by $\alpha(d\otimes d') = 1+d\smash d'$. Observe that $\End_R^*(F\smash F)$ is
complete with respect to the profinite filtration, because the $(F\smash F)_*$-module
\[
(F\smash F)_*^R(F\smash F) \cong (F\smash F)_*^R(F)\otimes_{(F\smash F)_*} (F\smash
F)_*^R(F)
\]
is free (compare \cite{jw}*{Remark 2.23}). Composition $\circ$ in
$\End_R^*(F\smash F)$ is clearly continuous, and so
$(\End_R^*(F\smash F), \circ)$ is a complete topological monoid.
It is easily checked that $\alpha$ is a continuous homomorphism of
topological  monoids. Moreover, the action of
$(\DDer^*_R(F)\otimes\DDer^*_R(F))^0$ on $\Prod_R(F)$ induced by
$\bar \pi$ is compatible with the canonical right action of
$\End_R^*(F\smash F)$ on $F^*_R(F\smash F)$ via $\alpha$. Since
$\End_R^*(F\smash F)$ is complete, $\alpha$ lifts to a continuous
homomorphism
\[
\wh\alpha\: (\DDer^*_R(F)\wh\otimes\DDer^*_R(F))^0 \to
\End_R^*(F\smash F).
\]

For the construction of $\hat\pi$, it remains to show that this
action restricts to an action on $\Prod_R(F)$. For this, we use
the facts that i) the action of $\End_R^*(F\smash F)$ on
$F^*_R(F\smash F)$ is continuous and ii) that $\Prod_R(F)$ is
closed in $F^*_R(F\smash F)$. Fact i) is easily verified. To prove
ii), we consider
\[
a\: F^*_R(F\smash F) \lra F^*_R(F\smash F\smash F), \ a(f) =
f(f\smash 1) - f(1\smash f),
\]
and the homomorphisms
\[
l, r\: F^*_R(F\smash F) \lra F^*_R(F), \ l(f)=f(1\smash \eta_F), \
r(f)=f(\eta_F\smash 1),
\]
where we implicitly use the equivalences $R\smash F\simeq F \simeq
F\smash R$. Observe that
\[
\Prod_R(F) = \ker(a)\cap\ker(l)\cap\ker(r)\cap F_R^0(F\smash
F)\subseteq F^*_R(F\smash F).
\]
Because $a, l$ and $r$ are continuous and because their targets
are Hausdorff, their kernels are closed. Moreover, so is
$F_R^0(F\smash F)$ and hence $\Prod_R(F)$.

It follows that $\Prod_R(F)$ is complete, as a closed subset of the complete module
$F^*_R(F\smash F)$. This implies that the action of
$(\DDer^*_R(F)\wh\otimes\DDer^*_R(F))^0$ on $F^*_R(F\smash F)$ restricts to an action on
$\Prod_R(F)$, and we are done.
\end{proof}

\section{Classification of products}\label{classification of prod}

In this section, we show hat the action of bilinear forms on
$I/I^2[1]$ on the set of products on a regular quotient $F=R/I$
classify the products on $F$. The main result is the following
theorem.

\begin{thm}\label{action}
Let $F=R/I$ be a regular quotient. Then the canonical action of
the group of bilinear forms $\Bil(I/I^2[1])$ on the set of
products $\Prod_R(F)$ is free and transitive.
\end{thm}

The strategy for the proof is as follows. On fixing
``coordinates'', we first give an explicit formula for $\beta\mu$,
for $\beta\in\Bil(I/I^2[1])$ and $\mu\in\Prod_R(F)$ (Lemma
\ref{explicitaction}). Secondly, we give an explicit description
of all products on $F$ (Lemma \ref{descriptionprod}). With these
two ingredients, we prove Theorem \ref{action}.

We first fix some notation. Let $F=R/I$ be a regular quotient and
let $\mu \in \Prod_R(F)$ be an arbitrary fixed product on $F$
(such a $\mu$ always exists, see \eg\ \cite{jw}*{Corollary 2.10}).
Let $(x_1, x_2, \ldots)$ be a regular sequence generating $I$.
Then the residue classes $\bar x_i\in V=I/I^2[1]$ form a basis,
and we let $\bar x_i^\vee\in D(V)$ denote the dual elements. An
arbitrary bilinear form $\beta \in \Bil(V)$ can be uniquely
written as a (possibly infinite) sum $\beta = \sum v_{ij}\, \bar
x_i^\vee \otimes \bar x_j^\vee$, with $v_{ij}=\beta(\bar x_i
\otimes \bar x_j) \in F_*$. Recall that $\psi\: \DDer^*_R(F) \to
D(V)$ from \eqref{defpsi} maps the Bockstein operation $Q_i$ to
$\bar x_i^\vee$, by definition of  $Q_i$.

Now $\bigl(\prod_{i+j\leq k} (1+ v_{ij}\, Q_i\smash Q_j)\bigr)_k$ is easily checked to be
a Cauchy sequence in the complete $F_*$-module $\End_R^*(F\smash F)$ (compare Section
\ref{def of action}). We define $\prod_{i,j} (1+v_{ij}\, Q_i\smash Q_j)$ to be its limit.

By definition of the canonical action of $\Bil(V)$ on
$\Prod_R(F)$, we have:

\begin{lem}\label{explicitaction}
In the notation from above, the product $\beta\mu$ is given by
\[
\beta \mu = \mu \circ \prod_{i,j} (1+v_{ij}\, Q_i \smash Q_j).
\]
\end{lem}

The next proposition describes the set of all products on $F$. It
is stated as Theorem 3.9 in \cite{ang}. We present a complete
proof here. It makes essential use of the existence of the
canonical action of $\Bil(I/I^2[1])$ on $\Prod_R(F)$, which in
turn relies crucially on the fact proved in \cite{jw} that
$\Der^*_R(F)$ is independent of the product on $F$, as a submodule
of $F^*_R(F)$.

\begin{lem}\label{descriptionprod}
For any product $\bar\mu \in \Prod_R(F)$, there exist uniquely determined elements
$v_{ij} \in F_*$ of degree $|v_{ij}|=|Q_i|+|Q_j|$ such that
\[
\bar\mu = \mu \circ \prod_{i,j}(1+v_{ij}\;  Q_i \smash Q_j).
\]
\end{lem}

The proof of Lemma \ref{descriptionprod} is postponed to the end
of the section. We first prove Theorem \ref{action}.

\begin{proof}[Proof of Theorem \ref{action}]
To prove transitivity, let $\mu, \bar\mu \in \Prod_R(F)$ be
arbitrary products. According to Lemma \ref{descriptionprod}, we
can write $\bar\mu$ as
\[
\bar\mu = \mu \circ \prod (1+v_{ij}\;  Q_i \smash Q_j).
\]
On setting $\beta=\sum v_{ij}\; \bar x_i^\vee \otimes \bar
x_j^\vee$, we obtain $\beta\mu = \bar\mu$, by Lemma
\ref{explicitaction}.

Freeness of the action follows from the fact that the coefficients
$v_{ij}$ in Lemma \ref{descriptionprod} are uniquely determined.
\end{proof}

We need some notation for the proof of Lemma
\ref{descriptionprod}. Let $a_i\in F_*^R(F)$ be the image of the
residue class $\bar x_i\in V$ under the characteristic
homomorphism $\phi\:V \to F_*^R(F)$. By \eqref{ffopexplicitintro},
we have $(F^\op)_*^R(F)\cong \Lambda (a_1,a_2,\ldots)$. Under this
isomorphism, $(F^\op)_*(Q_i)$ corresponds to the partial
derivative $\tfrac{\partial}{\partial a_i}$, see \cite{jw}*{Remark
4.5}. For a multi-index $I=(i_1,\ldots,i_m)$ with $i_1<\cdots
<i_m$, we write $|I|$ for $i_1+\cdots+i_m$, $Q_I$ for $Q_{i_1}
\cdots Q_{i_m}$ and $a_I$ for $a_{i_1}\wedge \cdots\wedge
a_{i_m}$.

\begin{proof}[Proof of Lemma \ref{descriptionprod}]
The K\"unneth homeomorphism (see \cite{sw}*{\S 2})
\[
\kappa\: F^*_R(F)\hatotimes F^*_R(F) \stackrel{\cong}{\lra}
F^*_R(F\wedge F)
\]
maps $\sum x_{IJ}\; Q_I \otimes Q_J$ to $\mu \circ(\sum x_{IJ}\;
Q_I \wedge Q_J)$. Since $\bar \mu \in F^0_R(F \wedge F)$, we may
write $\bar \mu =\mu \circ (\sum w_{IJ}\; Q_I \wedge Q_J)$, with
$|w_{IJ}|=|Q_I|+|Q_J|$. In particular, $w_{IJ}\neq 0$ only for
$|I|+|J|$ even. Since $\bar\mu$ has a two sided unit, it follows
that $w_{\emptyset J}=w_{I \emptyset}=0$ for all $I,J$. Hence
$\bar\mu$ can be written as
\[
\bar\mu = \mu\circ \biggl(1 +\sum_{|I|,|J|>0} w_{IJ}\; Q_I \wedge Q_J\biggr) =
\kappa\biggl(1 +\sum_{|I|,|J|>0} w_{IJ}\; Q_I \otimes Q_J\biggr).
\]

In a first step, we show that there exist $v_{IJ}\in F_*$ such
that
\begin{equation}\label{important}
1 +\sum_{|I|,|J|>0} w_{IJ}\; Q_I \wedge Q_J = \prod_{|I|,|J|>0}
\left(1+v_{IJ}\; Q_I \smash Q_J\right),
\end{equation}
where the product is taken in the monoid $\End_R^*(F\smash F)$.

If $(x_1,x_2,\ldots)$ is finite, the sum on the left hand side of \eqref{important} is of
the form $1 +\sum_{k=1}^{n} w_{I_{k}J_{k}}\; Q_{I_k} \wedge Q_{J_k}$. Set $\alpha_k
=w_{I_{k}J_{k}}\; Q_{I_k} \wedge Q_{J_k}$. By induction on $n$, one easily proves that
\begin{equation}\label{prodind}
\prod_{\substack{k=1,\ldots,n \\ 1 \leq i_1 <\cdots<i_k \leq n} }
(1+(-1)^{k-1}\alpha_{i_1}\cdots\alpha_{i_k})=1+\sum_{k=1}^n\alpha_k.
\end{equation}
This shows \eqref{important} for finite sequences $(x_1, x_2, \ldots)$. The general case
follows from this by passing to limits.

In a second step, we use the associativity of $\bar\mu$ to show that the coefficients
$v_{IJ}$ in \eqref{important} are zero for $|I|+|J|>2$. We write
\[
\bar\mu = \mu \circ \prod_{i,j} (1+v_{ij}\;  Q_i \smash Q_j)\circ
\prod_{|I|+|J|>2} (1+v_{IJ}\;  Q_I \smash Q_J)
\]
and let $\beta=\sum (-v_{ij})\; \bar x_i^\vee \otimes \bar
x_j^\vee \in \Bil(V)$. From Lemma \ref{descriptionprod}, we deduce
\[
\beta\bar\mu = \mu  \circ\! \prod_{|I|+|J|>2} (1+v_{IJ}\;  Q_I
\smash Q_J).
\]
We set $\tilde\mu=\beta\bar\mu \in \Prod_R(F)$ and assume that ${\mathcal I} = \{(I,J)\,
|\, v_{IJ}\neq 0\}$ is non-empty. We will show below that this implies that the two
morphisms
\begin{equation}\label{assoccond}
\tilde\mu_*(\tilde\mu_* \otimes 1), \;\tilde\mu_*(1 \otimes
\tilde\mu_*)\: (F^\op)_*^R(F)^{\otimes 3} \lra (F^\op)_*^R(F)
\end{equation}
are different, where $\tilde\mu_*$ stands for $m^{F^\op}_{\tilde
F}$. It follows that $\tilde\mu$ is not associative, which is a
contradiction. Therefore, $\mathcal I$ is empty, and the statement
is proved.

Let $(I_0,J_0) \in {\mathcal I}$ such that $|I_0|+|J_0|$ is
minimal. In the case where $|I_0|>1$, we set $I_0 = (L,M)$ with
$|L|,|M| \geq 1$. If $|I_0|=1$, we decompose $J_0$ in the same
way. We show that the two morphisms of \eqref{assoccond} don't
agree by evaluating them on $a_L \otimes a_M \otimes a_{J_0}$ if
$|I_0|>1$ or on $a_{I_0} \otimes a_L \otimes a_M$ if $|I_0|=1$.

As $(F^\op)_*^R(F) \cong \Lambda(a_1,a_2,\ldots)$, the set of elements $\{ a_I \otimes
a_J \otimes a_{K} \}_{I,J,K}$ forms a basis of the free $F_*$-module
$(F^\op)_*^R(F)^{\otimes 3}$. By minimality of $(I_0,J_0)$, we have $|I|>|L|$ or
$|J|>|M|$ for any $(I,J)\in\mathcal I$. This shows that
\begin{equation}\label{someobs}
F_*(Q_I) \otimes F_*(Q_J)(a_L \otimes a_M) = 0
\end{equation}
for all $(I,J) \in {\mathcal I}$. For $a,b\in F_*^R(F^\op)$, let
us write $a\smash b$ for $m^{F^\op}_F(a\otimes b)$. Using
\eqref{someobs}, we find:
\begin{align*}
\tilde\mu_*&(\tilde\mu_* \otimes 1)(a_L \otimes a_M \otimes a_{J_0}) = \tilde\mu_*(\mu_*
\otimes 1)(a_L \otimes a_M \otimes a_{J_0})
\\
& = \tilde\mu_*(a_L \wedge a_M \otimes a_{J_0}) = \tilde\mu_*(a_{I_0} \otimes a_{J_0}) =
\mu_*(a_{I_0}\otimes a_{J_0} - v_{I_0J_0} \cdot 1 \otimes 1)
\\
& = a_{(I_0,J_0)} - v_{I_0J_0} \cdot 1.
\end{align*}
Note that the negative sign appears because we let commute two
elements of odd degree. Similarly, we compute:
\begin{align*}
\tilde\mu_*&(1 \otimes \tilde\mu_*)(a_L \otimes a_M \otimes
a_{J_0})  = \tilde\mu_*(1 \otimes \mu_*)(a_L \otimes a_M
\otimes a_{J_0})\\
& = \tilde\mu_*(a_L \otimes a_M \wedge a_{J_0}) = \tilde\mu_*(a_{L} \otimes a_{(M,J_0)})
= \mu_*(a_{L}\otimes a_{(M,J_0)})  = a_{(I_0,J_0)}.
\end{align*}
This shows that the two morphisms in \eqref{assoccond} are different, as required.

Uniqueness of the coefficients $v_{ij}$ follows from the equality
\[
\tilde \mu_*(a_i \otimes a_j)= a_i \wedge a_j - v_{ij}\cdot 1,
\]
which we used in the argument above. This concludes the proof.
\end{proof}

\section{Transformation rules for the characteristic bilinear form}\label{action and
charact}

In this section, we describe how the action of the bilinear forms
affects characteristic bilinear forms and draw some consequences.

\begin{prop}\label{newbilinear}
Let $(F,k, \pi)$ be an admissible pair, where $F=R/I$ is a regular quotient ring. For a
bilinear form $\beta\in\Bil(I/I^2[1])$, we have
\[
b^k_{\beta F} = b^k_F - k_* \otimes \beta.
\]
\end{prop}

\begin{prop}\label{newbilinearII}
Let $(F,\bar F, 1_F)$ be an admissible pair, where $F, \bar F$ are regular quotient rings
with the same underlying quotient module $R/I$, endowed with two (possibly) different
products. For $\beta\in\Bil(I/I^2[1])$, we have
\[
 b^{\beta\bar F}_{F} = b^{\bar F}_{F} -\beta^t.
\]
\end{prop}

The proof of these two propositions  is technical and will be given at the end of this
section. We draw some consequences first.

\begin{cor}\label{newcharbil}
Let $b_F$ be the characteristic bilinear form of a regular
quotient ring $F=R/I$ and let $\beta\in\Bil(I/I^2[1])$ be a
bilinear form. Then the characteristic bilinear form of $\beta F$
is given by
\[
b_{\beta F} = b_F -(\beta + \beta^t).
\]
\end{cor}

\begin{proof}
The equalities of Propositions \ref{newbilinear} and
\ref{newbilinearII} imply that
\[
b_{\beta F}=b^{\beta F}_{\beta F}=b^{\beta F}_{F}-\beta=b^{F}_{F} -\beta^t - \beta =b_F -
(\beta + \beta^t).\qedhere
\]
\end{proof}

\begin{cor}\label{symetric}
The characteristic bilinear form  $b_F$ of a regular quotient ring $F$ is symmetric.
\end{cor}

\begin{proof}
Let $\mu$ denote the product on $F=R/I$. By \cite{jw}*{Corollary 2.10}, there exists a
diagonal product $\bar \mu$ on $F$ with respect to some regular sequence $(x_1, x_2,
\ldots)$ generating $I$. By Theorem \ref{action}, there exists $\beta \in\Bil(I/I^2[1])$
with $\beta \bar F = F$. Corollary \ref{newcharbil} implies that $b_F = b_{\bar F}
-(\beta + \beta^t).$ Now $b_{\bar F}$ is diagonal with respect to the basis $\bar x_1,
\bar x_2 , \ldots$ of $I/I^2[1]$ associated to the sequence $(x_1, x_2, \ldots)$
(\cite{jw}*{Prop.\@  2.35}). Therefore, $b_{F}$ is the sum of two symmetric bilinear
forms and therefore symmetric.
\end{proof}

\begin{cor}\label{relateffop}
For a regular quotient ring $F$ with characteristic bilinear form
$b_F$, we have $F^\op=b_F F$ and $b_{F^\op}=-b_F$. Therefore, $F$
is commutative if and only if $b_F=0$.
\end{cor}

\begin{proof}
As the bilinear forms $\Bil(I/I^2[1])$ act transitively on $\Prod_R(F)$, there exists
$\beta\in\Bil(I/I^2[1])$ with $ F^\op = \beta F$. Proposition \ref{newbilinear} implies
that $b_{F^\op}^F=b^F_{\beta F} = b_{F} - \beta$. But $b^F_{F^\op}$ is trivial by
\cite{jw}*{Prop.\@  2.21} and so $\beta=b_F$. From Corollary \ref{newcharbil}, we deduce
that $b_{F^\op} = b_F -(b_F + b_F^t)=-b_F$, since $b_F$ is symmetric.
\end{proof}

\begin{rem}\label{57}
Let $(F=R/I, \mu)$ be a regular quotient ring which is diagonal
with respect to some regular sequence $(x_1,x_2,\ldots)$
generating $I$. Then $b_F \in \Bil(I/I^2[1])$ is diagonal with
respect to the basis $\bar x_1, \bar x_2, \ldots$, as we used
above. Thus $b_F$ can be written as $\sum \alpha_i \bar x_i^\vee
\otimes \bar x_i^\vee$, where $\alpha_i \in F_*$ and where $\bar
x_i^\vee$ denotes the dual of $\bar x_i$. From Corollary
\ref{relateffop}, we obtain
\[
 \mu^{\op}= b_F\mu = \mu \circ \prod_{i} (1+\alpha_i\, Q_i \smash Q_i),
\]
where $Q_i$ denotes the Bockstein operation associated to $\bar x_i^\vee$. This
generalizes well-known formulas for $P(n)$ and $K(n)$ (see Section \ref{examples}).
\end{rem}

We now proceed to the proofs of Propositions \ref{newbilinear} and \ref{newbilinearII}.

Observe that it suffices to verify the statements for bilinear forms $\beta$ of the form
$\beta = \alpha\otimes\alpha'$ with $\alpha,\alpha'\in D(I/I^2[1])$, because an arbitrary
bilinear form can be written as a (possibly infinite) sum of bilinear forms of this type.

We first fix some notation used for the proofs. The proof of each proposition is then
preceded by a lemma.

Let $(F,k, \pi)$ be an admissible pair. For the proof of
Proposition \ref{newbilinearII}, $k$ will be $\bar F$ and
$\pi=1_F$. Let $\mu$ denote the product on $F$ and $\nu$ the one
on $k$. For $k=\bar F$, we write $\bar\mu$ instead of $\nu$, as
usual. We let $V=I/I^2[1]$ and consider
$\beta=\alpha\otimes\alpha'\in\Bil(V)$, where $\alpha, \alpha'\in
D(V)$. We let $d,d'\in\DDer^*_R(F)$ be the derivations
corresponding under $\psi\: \DDer^*_R(F) \cong D(V)$ to $\alpha,
\alpha'$, respectively. By definition of the action of $\Bil(V)$
on $\Prod_R(F)$, we have (using notation from Section \ref{def of
action})
\begin{equation}\label{e2}
\beta\mu =(\alpha\otimes\alpha')\mu=\mu_{d,d'}=\mu (1+d\smash d').
\end{equation}
We write $\bar x, \bar y$ for the residue classes of elements
$x,y\in I$ in both $V$ and in $k_*\otimes_{F_*} V$. Recall that
$b^{k}_{F}$ is defined as $b^{k}_{F}(\bar x\otimes \bar y)=\nu_*
k_*(\pi)(\phi(\bar x)\cdot\phi(\bar y))$, where $\phi$ is the
characteristic homomorphism $\phi^k_F\: V\to k_*^R(F)$ and where
$a\cdot b = m^{k}_{F}(a\otimes b)\in k_*^R(F)$ for $a,b\in
k_*^R(F)$ \eqref{multalgebra}.

\begin{lem}\label{actprod}
Let $(F,k, \pi)$ be an admissible pair, where $F=R/I$ is a regular
quotient ring. For $\beta$ a bilinear form in $\Bil(I/I^2[1])$ and
$x,y\in I$, we have:
\[
m^{k}_{\beta F}(\phi(\bar x)\otimes\phi(\bar y)) = \phi(\bar x)
\cdot \phi(\bar y)- \pi_*(\beta(\bar x\otimes \bar y)) \cdot 1.
\]
\end{lem}

\begin{proof}
Let $\beta = \alpha\otimes\alpha'$ with $\alpha,\alpha'\in D(V)$.
Recall the definition of $m^{k}_{\beta F}$:
\begin{equation}\label{e1}
m^{k}_{\beta F}(\phi(\bar x)\otimes\phi(\bar y)) = (\nu \wedge \beta\mu)_*  (1 \wedge
\tau \wedge 1)_*  \zeta (\phi(\bar x)\otimes\phi(\bar y)),
\end{equation}
where $\zeta\:k^R_*(F) \otimes k^R_*(F) \to (k\wedge F\wedge
k\wedge F)_*$ is the canonical map and $\tau $ the switch map
$\tau \: F \wedge k \to k\wedge F$. From the definition of
$\beta\mu$, we deduce that
\begin{align*}
& m^{k}_{\beta F}(\phi(\bar x)\otimes\phi(\bar y)) =((\nu \wedge
\mu +\nu \wedge (\mu\circ d\wedge d'))_*) (1 \wedge \tau \wedge
1)_* \zeta (\phi(\bar x)\otimes\phi(\bar y))
\\
& = \phi(\bar x)\cdot \phi(\bar y)-(\nu \wedge \mu)_*  (1\smash\tau\smash 1)_* \zeta
\bigl((1\smash d)_*(\phi(\bar
x))\otimes (1\smash d')_*(\phi(\bar y))\bigr) \\
& = \phi(\bar x)\cdot \phi(\bar y)-(\nu \wedge \mu)_* (1\smash\tau\smash 1)_* \zeta
\bigl(k_*^R(d)(\phi(\bar x))\otimes k_*^R(d')(\phi(\bar y))\bigr).
\end{align*}
By \cite{jw}*{Lemma 4.11}, we have that $k_*^R(d)(\phi(\bar x))=\alpha(\bar x)\cdot 1$
and $k_*^R(d')(\phi(\bar y))=\alpha'(\bar y)\cdot 1$, which implies the statement.
\end{proof}

\begin{proof}[Proof of Proposition \ref{newbilinear}]
Let $\beta = \alpha\otimes\alpha'$, $\alpha,\alpha'\in D(V)$.
Lemma \ref{actprod} implies:
\begin{align*}
b^{k}_{\beta F}(&\bar x\otimes \bar y)  =\psi_*\bigl(m^{k}_{\beta F}(\phi(\bar
x)\otimes\phi(\bar y))\bigr) =\psi_*\bigl(\phi(\bar x)\cdot \phi(\bar y)-\alpha(\bar
x)\alpha'(\bar y)\cdot 1\bigr)
\\ & = b^{k}_{F}(\bar x\otimes \bar y) - \alpha(\bar x)\alpha'(\bar y). \qedhere
\end{align*}
\end{proof}



\begin{lem}\label{lemma 2}
For $F,\bar F$ as in Proposition \ref{newbilinearII} and $\beta\in
\Bil(I/I^2[1])$, we have:
\[
m^{\beta \bar F}_{F}(\phi(\bar x)\otimes\phi(\bar y)) = \phi(\bar
x)\cdot\phi(\bar y) - \beta(\bar x\otimes \bar y)\cdot 1.
\]
\end{lem}

\begin{proof} 
Let $\beta = \alpha\otimes\alpha'$ with $\alpha,\alpha'\in D(V)$.
By definition of $\beta\bar\mu$, we have:
\begin{align*}\label{e3}
& m^{\beta \bar F}_{F}(\phi(\bar x)\otimes\phi(\bar y))
 = (\bar\mu \wedge \mu +(\bar\mu \circ d\wedge d') \wedge \mu)_* (1 \wedge \tau \wedge
1)_* \zeta (\phi(\bar x)\otimes\phi(\bar y))
\\
& = \phi(\bar x)\cdot \phi(\bar y)+(\bar\mu \wedge \mu)_* (1\smash\tau\smash 1)_* \zeta
\bigl((d\wedge 1)_* \otimes (d' \wedge 1)_*(\phi(\bar x)\otimes\phi(\bar y))\bigr)
\\
& = \phi(\bar x)\cdot \phi(\bar y)-(\bar\mu \wedge \mu)_*  (1\smash\tau\smash 1)_* \zeta
\bigl((d\wedge 1)_*(\phi(\bar x))\otimes(d' \wedge 1)_*(\phi(\bar y))\bigr).
\end{align*}
It remains to identify the second summand of the last equality
above. By definition, we have $(1\smash d)_*=F^R_*(d)$, and
furthermore
\[
(d\smash 1)_* = \tau_*(1\smash d)_*\tau_*=\tau_* F^R_*(d)\tau_*.
\]
>From \cite{jw}*{Prop.\@  3.6}, we obtain $\tau_*\phi(\bar x)=-\phi(\bar x)$. By
\cite{jw}*{Lemma 4.11}, we have $F_*^R(d)(\phi(\bar x))=\alpha(\bar x)\cdot 1$. This
yields $(d\smash 1)_* (\phi(\bar x))=-\alpha(\bar x)\cdot 1$. Analogously, we obtain
$(d'\smash 1)_* (\phi(\bar y))=-\alpha'(\bar y)\cdot 1$, and we are done.
\end{proof}

\begin{proof}[Proof of Proposition \ref{newbilinearII}] 
For $\beta = \alpha\otimes\alpha'$ with $\alpha,\alpha'\in D(V)$,
we compute:
\begin{align*}
b^{\beta\bar F}_{F}(&\bar x \otimes \bar y) =(\beta\bar\mu)_*(m^{\beta\bar F
}_{F}(\phi(\bar x)\otimes\phi(\bar y))) = (\beta\bar\mu)_*(\phi(\bar x)\cdot \phi(\bar
y)-\alpha(\bar x)\alpha'(\bar y)\cdot 1)
\\
& = (\bar\mu + \bar\mu (d\wedge d'))_*(\phi(\bar x)\cdot \phi(\bar
y))-\alpha(\bar x)
\alpha'(\bar y)\cdot 1)
\\
&= b^{\bar\mu}_{\mu}(\bar x \otimes \bar y)- \alpha(\bar x)\alpha'(\bar y)+(\bar\mu
(d\wedge d'))_*(\phi(\bar x)\cdot \phi(\bar y))).
\end{align*}
The first equality holds by definition of the characteristic
bilinear form, the second by Lemma \ref{actprod}, the third by
definition of $\beta\bar\mu$ and the fourth because $d$ and $d'$
are derivations and so are trivial on $1$.

Since $\alpha(\bar x)\alpha'(\bar y) = \beta(\bar x \otimes \bar y),$ it remains to show
that
\begin{equation*}
(\mu (d\wedge d'))_*(\phi(\bar x)\cdot \phi(\bar y)))=\beta(\bar x \otimes \bar
y)-\beta^t(\bar x \otimes \bar y).
\end{equation*}
To prove this, we write $d\smash d'$ as $(d\smash 1)(1\smash d')$.
Using computations from the proof of Lemma \ref{lemma 2} and the
fact that
$F^R_*(d)$ and $F^R_*(d')$ are derivations with respect to
$m^{\bar F}_{F}$ (see \cite{jw}*{Lemma 4.3}), we obtain:
\[
\begin{split}
(&\mu(d\smash d'))_*(\phi(\bar x) \cdot \phi(\bar y)) = \bar\mu_*(d\smash 1)_*(1\smash
d')_*(\phi(\bar x) \cdot \phi(\bar y))
\\
& = \bar\mu_*(d\smash 1)_* (\alpha'(\bar x)  \phi(\bar y) - \phi(\bar x)  \alpha'(\bar
y))  = \bar\mu_*(-\alpha'(\bar x)\alpha(\bar y) \cdot 1 + \alpha(\bar x)\alpha'(\bar
y)\cdot  1)
\\
& = (\alpha\otimes\alpha')(\bar x\otimes \bar y -\bar y\otimes \bar x) =\beta(\bar x
\otimes \bar y)-\beta^t(\bar x \otimes \bar y),
\end{split}
\]
and the proposition is proven.
\end{proof}

\section{Maps of quotient ring spectra}\label{naturality}

In this section, we determine which maps $\pi \: F \to G$ between
regular quotient rings are multiplicative. We start with a
definition.

\begin{defn}
An admissible pair $(F,G,\pi)$ with $F=R/I$ and $G=R/J$ is called
{\it smooth} if the canonical homomorphism $\pi_* \: G_*
\otimes_{F_*} I/I^2[1] \to J/J^2[1]$ is injective. If there is no
risk of confusion, we say that $I \subset J$ is smooth.
\end{defn}

\begin{thm}\label{nat}
Let $(F,G,\pi)$ be an admissible pair for which $F=R/I$ and $G=R/J$ are regular quotient
rings and which is smooth. Then $\pi$ is multiplicative if and only if $G_* \otimes
b_{F}=b^{G}_{F}=\pi^*(b_G)$.
\end{thm}

The strategy for the proof is as follows. We first prove the result in the special case
where $F$ is diagonal. As in this case the smoothness hypothesis is unnecessary, we
formulate a separate statement (Proposition \ref{natdiag}). After assembling some
auxiliary results (Lemmas \ref{technic} and \ref{technicii}), we prove Theorem \ref{nat}
by reducing it to the case where $F$ is diagonal.

\begin{prop}\label{natdiag}
Let $(F,G,\pi)$ be an admissible pair for which $F=R/I$ and $G=R/J$ are regular quotient
rings. Assume that $F$ is diagonal. Then $\pi$ is multiplicative if and only if $G_*
\otimes b_{F}=b^{G}_{F}=\pi^*(b_G)$.
\end{prop}

\begin{proof}
If $\pi$ is multiplicative, $(F,G,\pi)$ is a multiplicative
admissible pair by definition and the assertion follows from
\cite{jw}*{Prop.\@  2.20}.

To prove the converse, fix a regular sequence $(x_1,x_2,\ldots)$
generating $I$, for which there are products $\mu_k$ on the
$R/x_k$ such that the product $\mu_F$ on $F$ is the smash product
of the $\mu_k$. Let $\pi_k$ stand for the composition $\pi j_k \:
R/x_k \to F \to G$, where $j_k \: R/x_k \to F$ is the canonical
map.

By \cite{strickland}*{Prop.\@  4.8}, the map $\pi \: F \to G$ is multiplicative if and
only if i) \nolinebreak all the $\pi_k$ are multiplicative and ii) \nolinebreak $\pi_k$
commutes with $\pi_l$ for $k\neq l$.

In a first step, we show that the $\pi_k$ are multiplicative, \ie\
that they satisfy $\mu_G (\pi_k \smash \pi_k)=\pi_k \mu_k$, where
$\mu_G$ is the product on $G$. Because $x_k \in I \subset J$, the
$G_*$-module $G_*^R(R/x_k)$ is free on $1$ and
$a_{k}=\phi^G_{R/x_k}(\bar x_k)$, where $\phi^G_{R/x_k}$
is the characteristic homomorphism of the admissible pair $(R/x_k,
G, \pi_k)$. Therefore, the Kronecker duality homomorphism (see
\eg\ \cite{jw}*{Prop.\@ 2.25})
\[
d \: G^*_R(R/x_k\smash R/x_k) \lra \Hom^*_{G_*}(G_*^R(R/x_k\smash R/x_k), G_*)
\]
is an isomorphism. To relieve the notation, we identify $G_*^R(R/x_k\smash R/x_k)$ with
$G_*^R(R/x_k)\otimes G_*^R(R/x_k)$ via the K\"unneth isomorphism with respect to $\mu_G$
in the following.

To show that $\pi_k$ is multiplicative, we need to verify that
$d(\mu_G (\pi_k \smash \pi_k))$ and $d(\pi_k \mu_k)$ take the same
values on the basis elements $1 \otimes 1$, $1 \otimes a_{k}$,
$a_{k} \otimes 1$ and $a_{k} \otimes a_{k}$ of
$G_*^R(R/x_k)\otimes G_*^R(R/x_k)$. By naturality of the
characteristic homomorphism, we have
$G_*^R(\pi_k)(a_{k})=\phi_G(\bar\pi_*(\bar x_k))\in G_*^R(G)$,
where $\bar\pi_*\: G_*\otimes_{F_*}I/I^2[1]\to J/J^2[1]$ is
induced by $\pi_*$. Writing $a'_k$ for this element and
suppressing K\"unneth isomorphisms from the notation, we compute:
\begin{equation*}
\begin{split}
d(&\mu_G (\pi_k \smash \pi_k))(a_{k} \otimes a_{k}) = (\mu_G)_*G^R_*(\mu_G)(G^R_*(\pi_k)
\otimes G^R_*(\pi_k))(a_{k} \otimes a_{k})
\\
&=(\mu_G)_*G^R_*(\mu_G)(a'_{k} \otimes a'_{k}) = b_G(\pi_*(\bar x_k) \otimes \pi_*(\bar
x_k)) = \pi^*(b_G)(\bar x_k \otimes  \bar x_k).
\end{split}
\end{equation*}
On the other hand, we have (denoting both the residue classes of $x_k$ in $G_*
\otimes_{R_*/x_k} (x_k)/(x_k)^2[1]$ and in $G_* \otimes_{F_*} I/I^2[1]$ by $\bar x_k$):
\begin{equation*}
d(\pi_k \mu_k)(a_{k} \otimes a_{k}) = (\mu_G)_*G^R_*(\pi_k\mu_k)(a_{k} \otimes a_{k}) =
b_{R/x_k}^G(\bar x_k \otimes \bar x_k) = b^G_F(\bar x_k\otimes \bar x_k).
\end{equation*}
For the last equality, we have used that $j_k \: R/x_k \to F$ is multiplicative. By
hypothesis, we have $\pi^*(b_G) = b^G_F$, which shows that
\[
d(\mu_G (\pi_k \smash \pi_k))(a_{k} \otimes a_{k}) = d(\pi_k \mu_k)(a_{k} \otimes a_{k}).
\]
Similar, but simpler calculations show that $d(\mu_G (\pi_k \smash
\pi_k))$ and $d(\pi_k \mu_k)$ agree on the other basis elements $1
\otimes 1, 1 \otimes a_{k}$ and $a_{k} \otimes 1$ as well.

In a second step, we prove that $\pi_k$ and $\pi_l$ commute for $k
\neq l$, in the sense that $\mu_G (\pi_k\smash \pi_l) = \mu_G^\op
(\pi_k\smash\pi_l)$.
The relevant Kronecker duality morphism
\[
d\: G^*_R(R/x_k \smash R/x_l) \lra \Hom^*_{G_*}(G_*^R(R/x_k\smash R/x_l), G_*)
\]
is again an isomorphism. We use the notation and conventions from above and evaluate
$d(\mu_G^\op (\pi_k \smash \pi_l))$ and $d(\mu_G (\pi_k \smash \pi_l))$ on $a_{k} \otimes
a_{l}$. We first compute:
\begin{equation*}
\begin{split}
d(&\mu_G^\op (\pi_k \smash \pi_l))(a_{k} \otimes a_{l}) =
(\mu_G)_*G^R_*(\mu_G^\op)(G^R_*(\pi_k)\otimes G^R_*(\pi_l))(a_{k} \otimes a_{l})
\\
& = (\mu_G)_*G^R_*(\mu_G^\op)(a'_{k} \otimes a'_{l}) = (\mu_G)_*(a'_{k} \ast a'_{l}),
\end{split}
\end{equation*}
where $\ast$ denotes the product on $G_*^R(G^\op)$. Now $(\mu_G)_*
\: G_*^R(G^\op ) \to G_*$ is multiplicative by
\cite{jw}*{Corollary 3.3}. Together with $G_{\odd} =0$, this
implies that $(\mu_G)_*(a_{k} \ast a_{l}) = (\mu_G)_*(a_{k})\cdot
(\mu_G)_*(a_{l})=0$. On the other hand, we have:
\begin{equation*}
\begin{split}
d(& \mu_G (\pi_k \smash \pi_l))(a_{k} \otimes a_{l}) =(\mu_G)_*G^R_*(\mu_G)(a'_{k}
\otimes a'_{l}) =(\mu_G)_*(a'_{k} \cdot a'_{l})
\\
& = b_G(\bar\pi_*(\bar x_k) \otimes \bar\pi_*(\bar x_l))=\pi^*(b_G)(\bar x_k \otimes \bar
x_l),
\end{split}
\end{equation*}
where $\cdot$ denotes the product of $G_*^R(G)$. Since $\pi^*(b_G)=G_* \otimes b_F$ by
hypo\-thesis, since $b_F$ is diagonal with respect to the basis $\bar x_1, \bar
x_2,\ldots$ and since $k \neq l$, we have $\pi^*(b_G)(\bar x_k \otimes \bar x_l)=0$.

Leaving the analogous, simpler computations on $1 \otimes 1$, $1
\otimes a_{l}$, $a_{k} \otimes 1$ again to the reader, we conclude
that $d(\mu_G^\op (\pi_k \smash \pi_l))= d(\mu_G (\pi_k \smash
\pi_l))$. Hence $\pi_k$ and $\pi_l$ commute with each other, which
concludes the proof.
\end{proof}

By \cite{jw}*{Prop.\@  2.35}, the characteristic bilinear form of
a diagonal regular quotient ring is diagonal. We now show  that
the converse is true as well:

\begin{prop}\label{caradiag}
Let $F=R/I$ be a regular quotient ring and $(x_1,x_2,\ldots)$ a regular sequence
generating the ideal $I$. Then $F$ is diagonal with respect to the sequence
$(x_1,x_2,\ldots)$ if and only if $b_F$ is diagonal with respect to the basis $\bar
x_1,\bar x_2,\ldots$ of $I/I^2[1]$.
\end{prop}

\begin{proof}
The necessity of the condition was shown in \cite{jw}, as noted above. For sufficiency,
assume that $b_F$ is diagonal, and let $\mu_k$ be a product on $R/x_k$ such that the
canonical map $j_k \: R/x_k \to F$ is multiplicative, for all $k$. The proof of
Proposition \ref{natdiag} above shows that $j_k$ and $j_l$ commute if $k \neq l$, since
$b_F$ is diagonal with respect to the $\bar x_i$. From \cite{strickland}*{Prop.\@ 4.8},
we deduce that the product on $F$ is the smash ring product of the $\mu_k$.
\end{proof}

\begin{lem}\label{technic}
Let $(F,G,\pi)$ be an admissible pair satisfying the conditions of
Theorem \ref{nat}. Assume that $G_* \otimes b_{F}=\pi^*(b_G)$.
Then:

\begin{enumerate}\itemsep2pt

\item There exist products $\bar\mu$ on $F$ and
$\bar\nu$ on $G$ such that $\pi\: \bar F \to \bar G$ is
multiplicative.

\item For any $d \in \DDer^*_R(G)$ there exists  $\delta \in
\DDer^*_R(F)$ such that $d \pi = \pi \delta$.

\end{enumerate}
\end{lem}

\begin{proof}
(i) Let $\beta \in \Bil(I/I^2[1])$ be defined by $\beta(\bar x_i
\otimes \bar x_j) = 0$ for $i \geq j$, $\beta(\bar x_i \otimes
\bar x_j)= b_F(\bar x_i \otimes \bar x_j)$ for $i < j$ and let
$\tilde F = \beta F$. By Corollary \ref{newcharbil}, the
characteristic bilinear form $b_{\tilde F}$ of $\tilde F$ is given
by $b_{\tilde F}=b_F-(\beta + \beta^t)$ and is therefore diagonal
with respect to the $\bar x_i$.

Since $(F,G,\pi)$ is smooth, the homomorphism
\[
\pi^* \: \Bil(J/J^2[1]) \to \Bil(G_* \otimes_{F_*} I/I^2[1])
\]
is surjective. Choose $\gamma \in \Bil(J/J^2[1])$ with
$\pi^*(\gamma)= G_* \otimes \beta$ and set $\bar G= \gamma G$. By
hypothesis and by Corollary \ref{newcharbil}, it follows that $G_*
\otimes b_{\tilde F}=\pi^*(b_{\bar G})$.

Let $\pi_k= \pi j_k \: R/x_k \to \bar G$, with $j_k$ the canonical
map. The proof of Proposition \ref{natdiag} implies that $\pi_k$
and $\pi_l$ commute for $k \neq l$. Choose a product $\mu_k$ on
$R/x_k$ such that $\pi_k$ is multiplicative, for each $k$, and let
$\bar F$ be the induced smash ring spectrum. By
\cite{strickland}*{Prop.\@ 4.8}, $\pi \: \bar F \to \bar G$ is
then multiplicative.

(ii) Suppose first that $\pi$ is multiplicative. Then we have the
following commutative diagram:
\[
\begin{array}{c}
\xymatrix{\DDer^*_R(G)\ar[d]_-{\psi}^-{\cong} \ar[r]^-{-\circ \pi} &
\DDer^*_R(F,G)\ar[d]_-{\psi}^-{\cong}& \DDer^*_R(F) \ar[d]_-{\psi}^-{\cong}\ar[l]_-{\pi
\circ -}
\\
\Hom^*_{G_*}(J/J^2[1], G_*)  \ar[r]^-{\pi^*} & \Hom^*_{F_*}(I/I^2[1], G_*)
 & \Hom^*_{F_*}(I/I^2[1], F_*) \ar[l]_-{\pi_*} }
\end{array}
\]
where $\psi$ is as in \eqref{defpsi}. The right bottom morphism
$\pi_*$ is surjective, which implies the statement in this
particular case.

In the general case, (i) implies that there exist products $\bar
\mu$ on $F$ and $\bar \nu$ on $G$ such that $\pi \: \bar F \to
\bar G$ is multiplicative.  By \cite{jw}*{Lemma 4.6}, $d$ is a
derivation for {\em any} product on $G$, in particular $d \in
\DDer^*_R(\bar G)$. By what we have shown above, there exists
$\delta \in \DDer^*_R(\bar F)$ such that $d \pi = \pi \delta$. By
\cite{jw}*{Lemma 4.6} again, we deduce that $\delta \in
\DDer^*_R(F)$, which proves (ii).
\end{proof}

The following two statements are generalizations of Lemma
\ref{lemma 2} and Proposition \ref{newbilinearII}, respectively,
for the case where the map of the admissible pair is not
necessarily the identity.

\begin{lem}\label{technicii}
For an admissible pair $(F,G, \pi)$ satisfying the conditions of
Theorem \ref{nat} and $\gamma \in \Bil(J/J^2[1])$, we have
$b^{\gamma G}_F = b^G_F - \pi^*(\gamma^t).$
\end{lem}

\begin{proof}
Let $\phi = \phi^G_F$ be the characteristic homomorphism of the admissible pair $(F,G,
\pi)$. In a first step, we show that for $x,y\in I$, we have:
\begin{equation}\label{eqtech}
m^{\gamma G}_{F}(\phi(\bar x)\otimes\phi(\bar y)) = m^G_F(\phi(\bar x)\otimes\phi(\bar
y)) - \pi^*(\gamma)(\bar x\otimes \bar y)\cdot 1.
\end{equation}

Let $\mu$ be the product on $F$ and $\nu$ the one on $G$, and let
us write $a\cdot b$ for $m^{G}_{F}(a\otimes b) \in G^R_*(F)$,
where $a,b\in G_*^R(F)$. Clearly, it suffices to prove
\eqref{eqtech} for the case where $\gamma$ is of the form $\gamma
= \alpha \otimes \alpha'$, with $\alpha, \alpha' \in D(J/J^2[1])$.
Let $d,d' \in \DDer^*_R(G)$ correspond to $\alpha, \alpha'$,
respectively, under the isomorphism $\psi \: \DDer^*_R(G) \to
D(J/J^2[1])$. We have $\gamma \nu =\nu +\nu (d\wedge d')$. Recall
that for $x\in I$, we denote both the residue classes of $x\in I$
in $I/I^2[1]$ and in $G_* \otimes_{F_*} I/I^2[1]$ by $\bar x$.
Exactly as in the proof of Lemma \ref{lemma 2} we identify
$m^{\gamma G}_F(\phi(\bar x)\otimes\phi(\bar y))$ for $x,y\in I$
as
\begin{equation*}
\phi(\bar x)\cdot \phi(\bar y)-(\mu \wedge \nu)_* (1\smash\tau\smash 1)_* \zeta
\bigl((d\wedge 1)_*(\phi(\bar x))\otimes(d' \wedge 1)_*(\phi(\bar y))\bigr).
\end{equation*}
To determine $(d\wedge 1)_*(\phi(\bar x))$, we proceed as follows. By Lemma
\ref{technic}(ii), there exists $\delta \in \DDer^*_R(F)$ such that $\pi \delta =d \pi$.
By commutativity of the diagram
\begin{equation*}
\begin{array}{c}
\xymatrix{ F^R_*(F) \ar[d]_-{(\pi\wedge 1)_*}\ar[r]^-{(\delta
\wedge 1)_*} & F^R_*(F)\ar[d]^-{(\pi \wedge 1)_*}
\\
G^R_*(F) \ar[r]^-{(d \wedge 1)_*}  & G^R_*(F),  }
\end{array}
\end{equation*}
and using similar arguments as in the proof of Lemma \ref{lemma
2}, we deduce:
\begin{equation*}
\begin{split}
(d \wedge 1)_*(\phi(\bar x))& =(d \wedge 1)_*(\phi^G_F(\bar x))=(\pi \wedge 1)_*(\delta
\wedge 1)_*(\phi_F(\bar x))
\\
& = (\pi\wedge 1)_* \tau_*(1\wedge\delta)\tau_*(\phi_F(\bar x)) = -(\pi\wedge 1)_*
(\psi(\delta)(\bar x)\cdot 1)
\\
& = - \pi_*(\psi(\delta)(\bar x)) \cdot 1 = -\pi^*(\psi(d))(\bar
x) \cdot 1 = -\pi^*(\alpha)(\bar x)\cdot 1.
\end{split}
\end{equation*}
Similarly, we obtain $(d' \wedge 1)_*(\phi(\bar y))=
-\pi^*(\alpha')(\bar y)\cdot 1$. It follows that
\[
m^{\gamma G}_F(\phi(\bar x)\otimes\phi(\bar y))  = \phi(\bar
x)\cdot \phi(\bar y)- \pi^*(\alpha)(\bar x)\pi^*(\alpha')(\bar
y)\cdot 1,
\]
which is \eqref{eqtech} for $\gamma = \alpha\otimes\alpha'$.

We now proceed to the proof of the lemma itself. Again, we assume
$\gamma = \alpha\otimes\alpha'$, with $\alpha,\alpha'\in
D(J/J^2[1])$, and let $d=\psi(\alpha)$, $d'=\psi(\alpha')\in
\DDer^*_R(G)$. Using \eqref{eqtech}, we start the computation of
$b^{\gamma G}_{F}(\bar x \otimes \bar y)$ for $x, y\in I$ as in
the proof of Proposition \ref{newbilinearII} and find:
\begin{align*}
b^{\gamma G}_{F}(\bar x \otimes \bar y) &=b^{G}_{F}(\bar x \otimes \bar y)-
\pi^*(\gamma)(\bar x \otimes \bar y)+(\nu (d\wedge d')(1\smash\pi))_*(\phi(\bar x)\cdot
\phi(\bar y)).
\end{align*}
We now identify the last summand of the sum on the right hand
side. Since $d'\in \DDer^*_R(G)$ is a derivation, the homomorphism
\[
(1 \wedge d')_*=G^R_*(d') \: G_*^R(G) \to G_*^R(G)
\]
is a derivation, too \cite{jw}*{Lemma 4.3}. Using \cite{jw}*{Lemma
4.11} and writing $\cdot$ for $m^G_G$ (as well as for $m^G_F$), we
compute:
\begin{equation*}
\begin{split}
(\nu &(d\smash d')(1\smash\pi))_*(\phi(\bar x) \cdot \phi(\bar y))  = \nu_*(d\smash
1)_*(1\smash d')_*(\phi_G(\bar x) \cdot \phi_G(\bar y))
\\
& = \nu_*(d\smash 1)_*(G^R_*(d')(\phi_G(\bar x)) \cdot \phi_G(\bar y)-\phi_G(\bar x)
\cdot G^R_*(d')(\phi_G(\bar y)))
\\
& =\nu_*(d\smash 1)_*(\pi^*(\alpha')(\bar x)\phi_G(\bar y)-\phi_G(\bar
x)\pi^*(\alpha')(\bar y)).
\end{split}
\end{equation*}
In the proof of \eqref{eqtech} above, we showed that $-(d \wedge
1)_*(\phi(\bar x)) = -\pi^*(\alpha)(\bar x)\cdot 1$. Using the
analogous expression for $(d\wedge 1)_*(\phi(\bar y))$, we find
that
\[
b^{\beta G}_{F}(\bar x \otimes \bar y) = b^{G}_{F}(\bar x \otimes \bar y)-
\pi^*(\alpha')(\bar x) \pi^*(\alpha)(\bar y).
\]
This finishes the proof of the lemma.
\end{proof}

\begin{proof}[Proof of Theorem \ref{nat}]
If $\pi $ is multiplicative, then $(F,G,\pi)$ is a multiplicative admissible pair and the
statement follows from \cite{jw}*{Prop.\@  2.20}.

Conversely, let us assume that $b^{G}_{F}=\pi^*(b_G)=G_* \otimes b_{F}$. Let $\mu$ be the
product on $F$ and $\nu$ the one on $G$. Let $(x_1,x_2,\ldots)$ be a regular sequence
generating the ideal $I$. Let $\bar \mu$ be a product on $F$ which is diagonal with
respect to $(x_1, x_2, \ldots)$ (see \eg\ \cite{jw}*{Corollary 2.10}) and let $\beta \in
\Bil(I/I^2[1])$ be such that $\bar F= \beta F$. Write $\beta$ as a sum $\sum_i \epsilon_i
\otimes \epsilon_i'$ with $\epsilon_i, \epsilon_i' \in D_{F_*}(I/I^2[1])$. Because $(F,G,
\pi)$ is smooth, the composition
\[
\pi^* \:D_{G_*}(J/J^2[1]) \to D_{G_*}(G_* \otimes I/I^2[1])\cong
G_* \otimes_{F_*}D_{F_*}(I/I^2[1])
\]
is surjective. Choose $\alpha_i, \alpha_i' \in D_{G_*}(J/J^2[1])$ such that
$\pi^*(\alpha_i) = \epsilon_i$ and $\pi^*(\alpha_i')= \epsilon_i'$ and define $\gamma =
\sum \alpha_i \otimes \alpha_i'$. Observe that $\pi^*(\gamma)= G_* \otimes \beta$.

Now set $\bar G = \gamma G$. Using Proposition \ref{newbilinear}
and Lemma \ref{technicii}, we compute:
\[
b_{\bar F}^{\bar G}=b_{\beta F}^{\gamma G}=b_{ F}^{G}-
\pi^*(\gamma)^t-G_* \otimes \beta =\pi^*(b_G-\gamma^t - \gamma)=
\pi^*(b_{\bar G}).
\]
Similarly, we find $b_{\bar F}^{\bar G}=G_* \otimes b_{\bar F}$, and so
$G_* \otimes b_{\bar F}=b_{\bar F}^{\bar G}=\pi^*(b_{\bar G})$. Since $\bar F$ is
diagonal, this implies by Proposition \ref{natdiag} that $\pi \: \bar F \to \bar G$ is
multiplicative.

Let $d_i, d_i'\in\DDer^*_R(G)$ be the derivations corresponding to
$\alpha_i,\alpha_i'$ under $\psi\:\DDer^*_R(G) \cong
D_{G_*}(J/J^2[1])$, and $\delta_i, \delta_i'\in\DDer^*_R(F)$ under
$\psi\: \DDer^*_R(F) \cong D_{F_*}(I/I^2[1])$ to $\epsilon_i,
\epsilon'_i$. By naturality of $\psi$, we have $d_i \pi=\pi
\delta_i$ and $d_i' \pi=\pi \delta_i'$. From the definition of the
canonical action of the group of bilinear forms on the set of
products, we have that
\begin{align*}
\gamma \bar \nu \circ (\pi \wedge \pi)&=\bar\nu \circ  \prod_i(1 +
d_i \wedge d_i') (\pi \wedge \pi) =\bar\nu  \circ (\pi \wedge \pi)
\circ  \prod_i(1 + \delta_i \wedge
\delta_i')\\
&=\pi  \bar \mu\circ \prod_i(1 + \delta_i \wedge \delta_i')=\pi \circ \beta \bar\mu.
\end{align*}
Therefore $\pi \: F= - \beta \bar F \to - \gamma \bar G = G$ is
multiplicative. This completes the proof of the theorem.
\end{proof}

\section{Classification of products up to equivalence}\label{prod up to iso}

In this section, we classify the products on regular quotients up to equivalence.
Moreover, we study commutative products and consider the question of diagonalizability of
products on regular quotients.

Let $F=R/I$ be a regular quotient ring with product $\mu$. If
$\bar\mu$ is a second product on $F$, we write $\bar F$ for $F$,
endowed with $\bar \mu$, as before. If $\beta \in \Bil(I/I^2[1])$
is such that $\bar\mu=\beta \mu$, we alternatively write $\bar F =
\beta F$.

Recall the following definition:

\begin{defn}
Two products $\mu$ and $\bar \mu$ on $F$ are {\it equivalent} (denoted $\mu \sim \bar
\mu$) if there is a multiplicative isomorphism $f \: F \to \bar F$ in $\DR$. Such a map
$f$ is called a {\it multiplicative equivalence.}
\end{defn}

Together with Theorem \ref{action}, the following result gives a
classification for products up to equivalence:

\begin{thm}\label{altiso}
Let $F=R/I$ be a regular quotient ring and
$\beta\in\Bil(I/I^2[1])$ a bilinear form. Then $F$ and $\beta F$
are equivalent if and only if $\beta$ is alternating. In this
case, there is a canonical multiplicative equivalence $F\to \beta
F$.
\end{thm}

Let $F=R/I$ be a regular quotient ring. Consider the map
\[
\theta \: (\DDer^*_R(F)\times\DDer^*_R(F))^0 \to F^0_R(F)
\]
defined by $\theta(d,d')=1+dd'$. Since $F^*_R(F) \cong
\wh\Lambda(\DDer^*_R(F, F))$ (by \eqref{cohomologyiso}), the image
of $\theta$ is contained in the center of the monoid $F^*_R(F)$,
the product on $F^*_R(F)$ being the composition. Clearly, $\theta$
is bilinear. Since $F^*_R(F)$ is complete with respect to the
profinite topology, $\theta$ induces (see \eqref{bilder})
\[
\Theta\: \Bil(I/I^2[1]) \cong
(\DDer^*_R(F)\wh\otimes\DDer^*_R(F))^0 \lra F^0_R(F).
\]
The next lemma is a crucial step in the proof of Theorem
\ref{altiso}.


\begin{lem}\label{bilalt}
Let $F=R/I$ be a regular quotient ring and $\beta \in
\Alt(I/I^2[1])$. Then $\Theta(\beta)$ is a multiplicative
equivalence $\Theta(\beta)\: F\to \beta F$.
\end{lem}

\begin{proof}
It suffices to prove the lemma for bilinear forms $\beta$ of the form $\beta =
\alpha\otimes\alpha'-\alpha'\otimes\alpha$ with $\alpha,\alpha'\in D(I/I^2[1])$, because
an arbitrary alternating bilinear form can be written as a sum of such elements. Let
$d,d' \in \DDer^*_R(F)$ correspond to $\alpha, \alpha'$ under the isomorphism $\psi\:
\DDer_R^*(F) \cong D(I/I^2[1])$ \eqref{defpsi}. Denoting by $\mu$ the product on $F$, we
then have $\beta\mu=\mu (1+ d\wedge d')(1-d'\wedge d) $. In order to simplify the
notation, we write $\bar\mu$ for $\beta \mu$.

Since derivations anticommute, the map $f = 1 + d d'$ is an
equivalence, with inverse $1-dd'$. We have to show that $f \: F
\to  \bar F$ is multiplicative, that is, $f \mu = \bar \mu ( f
\smash f)$. For this, we first compute:
\[
\begin{split}
f \mu &=(1+ dd') \mu = \mu\(1+(d \wedge 1 + 1 \wedge d)(d' \wedge
1 + 1 \wedge d')   \)
\\
&= \mu\(1+dd' \wedge 1 + d \wedge d' - d' \wedge d +1 \wedge dd'\).
\end{split}
\]
On the other hand, we find:
\[
\begin{split}
\bar \mu ( f & \smash f)=\mu (1+ d\wedge d')(1-d'\wedge d) (1+dd' \wedge 1 + 1 \wedge
dd'+ dd' \wedge dd'   )
\\
& = \mu(1+dd' \wedge 1 +1 \wedge dd'  +dd' \wedge dd'  - d' \wedge d + d \wedge d'+dd'
\wedge d'd ) .
\end{split}
\]
Since $dd'=-d'd$, the lemma is proven.
\end{proof}

\begin{proof}[Proof of Theorem \ref{altiso}]
We fix a regular sequence $(x_1,x_2,\ldots)$ generating $I$. The
residue classes $\bar x_1,\bar x_2,\ldots$ form a basis of
$V=I/I^2[1]$, and we denote by $\bar x_1^\vee,\bar
x_2^\vee,\ldots$ the elements dual to the $\bar x_i$. The
Bockstein operations $Q_i$ are defined as $Q_i = \psi^{-1}(\bar
x_i^\vee)$, where $\psi$ is the isomorphism $\psi\: \DDer^*_R(F)
\to D(V).$


Assume first that $\beta$ is alternating. Then it can be written
as $\beta = \sum v_{ij} \bar x_i^\vee \otimes \bar x_j^\vee$ with
$v_{ii}=0$ and $v_{ij}=-v_{ji}$ for $i \neq j$. As a consequence,
$\beta \mu$ can be expressed as (see Lemma \ref{explicitaction}):
\begin{equation}\label{prodalt}
\beta \mu= \mu \prod_{i<j} \bigl((1+v_{ij} Q_i \smash Q_j) (1 - v_{ij}Q_j \smash Q_i)
\bigr).
\end{equation}

If the product in \eqref{prodalt} is finite, the map
$f=\prod_{i<j} (1 + v_{ij} Q_iQ_j)$ is a multiplicative homotopy
equivalence $f\: F\to \beta F$ by Lemma \ref{bilalt}. If the
product in \eqref{prodalt} is infinite, the Cauchy sequence of
multiplicative equivalences $\bigl(\prod_{i<j,\, i+j\leq k} (1 +
v_{ij} Q_iQ_j)\bigr)_k$ converges to one from $F$ to $\beta F$.

Suppose now that $F$ and $\bar F =\beta F$ are equivalent via a
multiplicative equivalence $\pi \: F \to\bar F$. Since $\pi_* \:
F_*\to \bar F_*$ and the induced homomorphism $\bar\pi_* \: \bar
F_* \otimes I/I^2[1] \cong I/I^2[1] \to I/I^2[1]$ are (equivalent
to) the identities, naturality of the characteristic bilinear form
and the commutative diagram
\[
\begin{array}{c}
\xymatrix{ F \ar[d]_-{\pi} \ar[r]^-{1_F}& F \ar[d]^-{\pi}
\\
\bar F  \ar[r]^-{1_{\bar F}}& \bar F }
\end{array}
\]
show that $b_F =b_{\bar F}$. Corollary \ref{newcharbil} implies that $b_{\bar F}= b_F
-(\beta + \beta^t)$. It follows that $\beta$ is antisymmetric. Hence it remains to check
that $\beta(\bar x_i \otimes \bar x_i)=0$ for all $i$ in order to prove that $\beta$ is
alternating.

Choose a product on $R/x_i$ such that the natural map $j_{i} \: R/x_i \to F$ is
multiplicative. Then both $(R/x_i,F, j_{i})$ and $(R/x_i,\bar F, \pi j_{i})$ are
multiplicative admissible pairs, and \cite{jw}*{Prop.\@  2.21} implies that
\[
b^F_{(R/x_i)^\op}=0=b^{\bar F}_{(R/x_i)^\op}.
\]
Since $(x_i) \subset I$ is smooth, Lemma \ref{technicii} applies.
On setting $b_i =b_{R/x_i}$ and recalling that $(R/x_i)^\op = b_i
R/x_i$ (Corollary \ref{relateffop}), we obtain:
\[
\begin{split}
0 & = b^{\bar F}_{(R/x_i)^\op}  = b^{\beta F}_{b_i R/x_i} = b^{\beta F}_{ R/x_i} - F_*
\otimes b_i = b^{F}_{ R/x_i} - F_* \otimes  b_i - j_{i}^*( \beta^t)
\\
& =   b^{F}_{b_i R/x_i} + j_{i}^*( \beta) = b^{F}_{(R/x_i)^\op} + j_i^*(\beta) =
j_{i}^*(\beta).
\end{split}
\]
Therefore $0=j_{i}^*( \beta)(\bar x_i \otimes \bar x_i)=\beta(\bar
x_i \otimes \bar x_i)$, where $\bar x_i$ again stands for the
residue class of $x_i$ in either $(x_i)/(x_i)^2[1]$ or $I/I^2[1]$.
Thus $\beta $ is alternating, and the theorem is proven.
\end{proof}

\begin{rem}\label{remclassprod}
Theorem \ref{altiso} states that $\Alt(I/I^2[1])$ acts freely and
transitively on the equivalence class of any product on $F$.
Therefore, the (additive) group of quadratic forms
$\Quad(I/I^2[1]) \cong \Bil(I/I^2[1]) / \Alt(I/I^2[1])$ acts
freely and transitively on the set of equivalence classes of
products on $F$.
\end{rem}

\begin{cor}\label{isoquad}
Let $\mu$ and $\bar\mu$ be two products on a regular quotient $F$.
\begin{enumerate}\itemsep2pt
\item If $F$ and $\bar F$ are equivalent then $b_F = b_{\bar F}$.
\item If $F_*$ is $2$-torsion-free, then  $F$ and $\bar F$ are
(canonically) equivalent if and only if $b_F = b_{\bar F}$ if and
only if $q_F = q_{\bar F}$.
\end{enumerate}
\end{cor}

\begin{proof}
(i) This has been shown in the proof of Theorem \ref{altiso}.

(ii) Suppose that $b_F = b_{\bar F}$. Let $\beta \in
\Bil(I/I^2[1])$ be the bilinear form which satifies $\beta F =
\bar F$ (Theorem \ref{action}). As in the proof of Theorem
\ref{altiso}, we deduce that $\beta^t =-\beta$. As $F_*$ is
$2$-torsion-free, this means that $\beta$ is alternating. Now
Theorem \ref{altiso} implies that $F$ and $\beta F=\bar F$ are
equivalent. The last implication is clear.
\end{proof}

\begin{rem}
If $F_*$ has $2$-torsion, there may exist non-equivalent products $\mu, \bar\mu \in
\Prod_R(F)$ with $b_F = b_{\bar F}$, see for instance Proposition \ref{classicalmorava}.
\end{rem}

\begin{rem}\label{bmap}
Let $F=R/I$ be a regular quotient ring. We may interpret the
characteristic bilinear form as a map $b \: \Prod_R(F) \to
\Bil(I/I^2[1])$. By Corollary \ref{symetric}, the image of $b$ is
contained in $\Sym(I/I^2[1]) \subset \Bil(I/I^2[1])$, and by
Corollary \ref{isoquad}, $b$ factors through the set of
equivalence classes of products on $F$:
\[
\bar b \: \Prod_R(F) / \sim \ \lra \Sym(I/I^2[1]) .
\]
>From Corollary \ref{isoquad} we deduce that $\bar b$ is injective if $F_*$ is
$2$-torsion-free. Moreover, we easily check that $\bar b$ is surjective if $2 \in F_*$ is
invertible.
\end{rem}

We now turn to a discussion of commutative products on a regular
quotient $F$. We prove that if $2\in F_*$ is invertible, there are
many commutative products in general (Proposition \ref{comprod}),
which, however, are all equivalent to each other (Corollary
\ref{commprodequiv}).

Strickland proves \cite{strickland}*{Theorem 2.6} that $F_*$ is
strongly realizable \cite{strickland}*{Def.\@ 2.1} if $2\in F_*$
is invertible. In particular, this shows that $F$ admits a
commutative product. Because $F_*$ is a quotient of $R_*$, any
commutative product on $F$ which is equivalent to a strong
realization is itself a strong realization. As a consequence, any
commutative product turns $F$ into a strong realization of $F_*$.

\begin{prop}\label{comprod}
Let $F=R/I$ be a regular quotient of $R$.

\begin{enumerate}\itemsep2pt
\item
Suppose that $F$ admits a commutative product. Then 
$\Asym(I/I^2[1])$ acts freely and transitively on the set of all
commutative products. \item If $2$ is invertible in $F_*$, there
exists a commutative product on $F$.
\end{enumerate}
\end{prop}

\begin{proof}
(i) Endow $F$ with a commutative product. For a bilinear form
$\beta \in \Bil(I/I^2[1])$, Corollaries \ref{newcharbil} and
\ref{relateffop} imply that $\beta F$ is commutative if and only
if $\beta \in \Asym(I/I^2[1])$. Now the statement follows from
Theorem \ref{action}.

(ii) Let $\mu$ be an arbitrary product on $F$ (see \eg\ \cite{jw}*{Corollary 2.10}) and
let $\beta = \tfrac{1}{2}\, b_F$. Then Corollary \ref{newcharbil} implies that
\[
b_{\beta F}= b_F -(\tfrac{1}{2}\, b_F + \tfrac{1}{2}\, b_F^t)=0,
\]
since $b_F$ is symmetric. Therefore $\beta F$ is commutative, by Corollary
\ref{relateffop}.
\end{proof}

\begin{rem}
Proposition \ref{comprod} is a generalization of \cite{strickland}*{Corollary 3.12},
which treats the case $F=R/x$. Note that in this situation, $\Asym(I/I^2[1])$ is the
group of 2-torsion elements in $F_{2|x|+2}$.
\end{rem}

Using Theorem \ref{altiso}, we deduce the following corollary.

\begin{cor}\label{commprodequiv}
Let $F=R/I$ be a regular quotient of $R$.
\begin{enumerate}\itemsep2pt
\item Suppose that $F$ admits a commutative product. Then the
group \linebreak $\Asym(I/I^2[1])/\Alt(I/I^2[1])$ acts freely and
transitively on the set of equivalence classes of commutative
products on $F$. \item If $F_*$ has no $2$-torsion, then there
exists at most one commutative product up to canonical
equivalence. \item If $2$ is invertible in $F_*$, there exists a
unique commutative product up to canonical equivalence.
\end{enumerate}
\end{cor}

\begin{rem}
If $F_*$ has 2-torsion, there may not exist any commutative product on $F$. This is
well-known, see \eg\ Proposition \ref{classicalmorava}.
\end{rem}

For regular quotients whose coefficient ring is $2$-torsion, we
have the following result.

\begin{prop}\label{712}
Let $F=R/I$ be a regular quotient such that $2 \cdot F_*=0$. Then
there exists a commutative product on $F$ if and only if $F$
admits a product whose characteristic bilinear form is
alternating. If this holds, then $b(\Prod_R(F)) = \Alt(I/I^2[1])$,
where $b \: \Prod_R(F) \to \Sym(I/I^2[1])$ is the map from Remark
\ref{bmap}.
\end{prop}
\begin{proof}
Assume that $F$ is endowed with a commutative product. For any
$\beta \in \Bil(I/I^2[1])$, Corollary \ref{newcharbil} implies
that $b_{\beta F}=b_F+\beta + \beta^t$. Hence  for any $\bar x \in
I/I^2[1]$ we have $b_{\beta F}(\bar x \otimes \bar x)=b_F (\bar x
\otimes \bar x) +2 \beta(\bar x \otimes \bar x)=0 $ since $b_F
=0$. As a consequence $b_{\beta F}\in\Alt(I/I^2[1])$ and thus
$b(\Prod_R(F)) \subset \Alt(I/I^2[1])$.

Conversely, let $F$ be endowed with a product such that $b_F \in
\Alt(I/I^2[1])$. Choose a regular sequence $(x_1, x_2,\ldots)$
generating the ideal $I$. Define $\beta \in \Bil(I/I^2[1])$ by
$\beta(\bar x_i \otimes \bar x_j)=0$ for $i \leq j$ and
$\beta(\bar x_i \otimes \bar x_j)=b_F(\bar x_i \otimes \bar x_j)$
for $i > j$. Then $b_{\beta F}=b_F+\beta + \beta^t$ is diagonal
with respect to the basis consisting of $\bar x_1, \bar
x_2,\ldots$. Since $b_F$ is alternating, this implies that
$b_{\beta F}(\bar x_i \otimes \bar x_i)=b_F (\bar x_i \otimes \bar
x_i)=0$ and hence that $b_{\beta F}=0$. From Corollary
\ref{relateffop}, it follows that $\beta F$ is commutative.

For the remaining statement, let $F$ be endowed with a commutative product and let $\beta
\in \Alt(I/I^2[1])$ be any alternating bilinear form. With the notation from above, we
define $\gamma \in \Bil(I/I^2[1])$ by $\gamma(\bar x_i \otimes \bar x_j)=0$ for $i \leq
j$ and $\gamma(\bar x_i \otimes \bar x_j)=\beta(\bar x_i \otimes \bar x_j)$ for $i > j$.
Then the characteristic bilinear form of $\gamma F$ satisfies $b_{\gamma F}= b_F +\gamma
+ \gamma^t=\beta$, and the proof is complete.
\end{proof}

We close this section with a discussion of diagonalizability of products. Recall
(\cite{jw}*{Def.\@ 2.9}) that a regular quotient ring $F$ is {\it diagonalizable} if it
is equivalent to a diagonal regular quotient ring.

Recall that the maximal ideal of a regular local ring of dimension
$n<\infty$ is always generated by a regular sequence of length $n$
(\cite{serre}*{Chap.\@ IV}).

\begin{prop}\label{diagprod}
Assume that $R_*$ is a regular local ring of dimension $n$ with
maximal ideal $I$ whose residue field $R_*/I$ is of characteristic
$p\geq 0$. Let $F=R/I$.

\begin{enumerate}

\item If $p$ is zero or an odd prime, then  $F$ is diagonalizable.

\item If $p=2$, then:
\begin{enumerate}
\item If $b_F \not\in \Alt(I/I^2[1])$, $F$ is  diagonalizable.
\item If $b_F \in \Alt(I/I^2[1])$ and $b_F\neq 0$, $F$ is not
diagonalizable. \item If $b_F =0$, $F$ is diagonalizable.
\end{enumerate}
\end{enumerate}
\end{prop}

\begin{proof}
Suppose first that (i) $p$ is zero or odd or that (ii) $p=2$ and $b_F \not\in
\Alt(I/I^2[1])$. Then \cite{bourbaki}*{Chap.\@ IX, § 6, Theorem 1} implies that there
exists a basis $\mathscr B$ consisting of elements $ b_0,\ldots,b_{n-1}$ of $I/I^2[1]$
such that the matrix of $b_{F}$ with respect to $\mathscr B$ is diagonal. By
\cite{serre}*{Chap.\@ IV, Prop.\@ 22}, there exists a regular sequence
$(y_0,\ldots,y_{n-1})$ generating $I$ such that $b_i=\bar y_i \in I/I^2[1]$ for all $i$.
We then conclude with Proposition \ref{caradiag}.

Now suppose that $p=2$ and $0 \neq b_F \in \Alt(I/I^2[1])$. For any basis $\mathscr B$
consisting of elements $b_0,\ldots,b_{n-1}$ of $I/I^2[1]$, we have $b_F(b_i \otimes
b_i)=0$ for all $i$. Therefore, $b_F$ is not diagonalizable, since $b_F \neq 0$. Hence,
by Proposition \ref{caradiag} again, $F$ is not diagonalizable.

For the remaining case, $p=2$ and $b_F=0$, the statement follows from Proposition
\ref{caradiag}. Alternatively, we may observe that $F$ is commutative (Corollary
\ref{relateffop}) and therefore diagonalizable (\cite{jw}*{Corollary 2.12}).
\end{proof}

\section{Examples}\label{examples}

In this section, we present some applications of our results.

We first collect some facts concerning complex cobordism. Let $MU$ be the commutative
$\SS$-algebra associated to complex cobordism (see \cite{ekmm}). Recall that there is an
isomorphism
\[
MU_* \cong \Z[x_1,x_2,\ldots],\; |x_i|=2i.
\]
Fix a prime number $p$ and recall from  \cite{strickland} that
$w_k\in MU_{2(p^k-1)}$ denotes the bordism class of a smooth
hypersurface $W_{p^k}$ of degree $p$ in $\C P^{p^k}$. Let
$J_n\subseteq MU_*$ be the ideal $(w_0, \ldots, w_{n-1})$, where
$w_0=p$. The following statement is an important ingredient for
our examples. It is a consequence of \cite{strickland}*{§7} and
\cite{jw}*{Prop.\@ 2.27}.

\begin{prop}\label{compofc}
Let $p=2$. There is a product on $F= MU/w_k$ with $b_F(\bar w_k
\otimes \bar w_k)\equiv w_{k+1} \mod J_k$ for $k\geq 0$.
\end{prop}

\subsection{$BP$-theory}

We fix a prime number $p$. The Brown-Peterson spectrum $BP$  can
be described as a regular quotient $BP=MU_{(p)}/I$ of the
$p$-localization $MU_{(p)}$ of $MU$, where $I \subset
(MU_{(p)})_*$ is the ideal generated by the regular sequence
$x_i,\; i \neq p^k-1 , \; k > 0$ (see \cite{ekmm} or
\cite{strickland}). The coefficient ring is given by
\[
BP_* \cong \Z_{(p)}[v_1,v_2,\ldots], \; |v_i|=2(p^i-1),
\]
where we choose the $v_i$'s to be Hazewinkel's generators (see
\cite{strickland}).

\begin{rem}\label{BPlocal}
Since $BP_*$ is $p$-local, we do not need to distinguish between
$MU$-products and $MU_{(p)}$-products on $BP$, see
\cite{ekmm}*{Section VIII. 3}.
\end{rem}

It is shown in \cite{strickland} that $BP$ is a commutative
$MU$-ring.

\begin{prop}\label{bp prod}
There are infinitely many non-equivalent $MU$-products on $BP$,
all of which induce the same ring structure in $\mathscr D_{\SS}$.
Infinitely many of the $MU$-products on $BP$ are commutative, but
all of these are equivalent.

\end{prop}

\begin{proof}
The equivalence classes of $MU$-products on $BP$ are in one-to-one
correspondence with the quadratic forms on $I/I^2[1]$. There are
infinitely many such, for odd $p$ for instance the ones associated
to the family of bilinear forms $\beta_k= v_k \, \bar
x_{i(k)}^\vee \otimes \bar x_{i(k)}^\vee$, where
$i(k)=\tfrac{1}{2}(p^k-1)$.

Let $\mu_0$ be a commutative product on $BP$ \cite{strickland}.
Any other product $\mu$ is of the form $\mu = \mu_0 \circ \prod(1+
a_{ij} Q'_i \wedge Q'_j)$, where $Q'_k\in BP^*_{MU_{(p)}}(BP)$ is
the Bockstein operation associated to $\bar x_k^\vee\in
D(I/I^2[1])$ (the notation $Q_k$ is reserved for a different
Bockstein operation, see the next section). Since $BP^*(BP)$ is
concentrated in even dimensions, all the $Q'_k$ are in the kernel
of the forgetful morphism
\[
BP^*_{MU_{(p)}}(BP) \lra BP^*(BP)
\]
induced by the monoidal functor $\mathscr{D}_{MU_{(p)}} \subset
\mathscr{D}_{\SS}$. As a consequence, all the $MU$-products on
$BP$ are equal to $\mu_0$ in $\mathscr{D}_{\SS}$.

The last assertion follows from Corollary \ref{commprodequiv}.
\end{proof}

\subsection{$P(n)$-theory}
We fix a prime number $p$ and endow $BP$ with a  commutative
$MU$-product. Recall that  $J_n\subseteq MU_*$ is the ideal $(w_0,
\ldots, w_{n-1})$, where $w_0=p$. The sequence of the $w_i$ is
regular, and the image of $J_n$ in $BP_*$ is the ideal $I_n=(v_0,
\ldots, v_{n-1})$, with $v_0=p$ (see \cite{strickland}).

We define $P(n)$ as a quotient of $BP$ (see \cite{strickland}):
\[
P(n) = BP/I_n = BP \wedge_{MU} MU / J_n.
\]
The coefficient ring satisfies $P(n)_* \cong \mathbb
F_p[v_n,v_{n+1},\ldots]$. The kernel $H_n$ of the composition
$(MU_{(p)})_* \to BP_* \to P(n)_*$ is generated by a regular
sequence. Therefore, $P(n)=MU_{(p)}/H_n$ is a regular quotient of
$MU_{(p)}$.

Since $P(n)_*$ is $p$-local, we do not need to distinguish between $MU$-pro\-ducts and
$MU_{(p)}$-products on $P(n)$ (see Remark \ref{BPlocal}).

We endow $P(n)$ with an $MU$-product $\mu_n$ as follows. If $p$ is
odd, $MU / J_n$ carries a commutative product $\nu$, since $2$ is
invertible. If $p=2$, we define a product $\nu$ on $MU / J_n$ as
the smash ring product of the $\nu_k$ of Proposition \ref{compofc}
for $k=0,\ldots ,n-1$. In any case, we define $\mu_n$ as the smash
ring product of $\mu_0$, a commutative product on $BP$, with
$\nu$. Observe that the natural map $\pi_n\: BP \to P(n)$ is then
multiplicative.

\begin{prop}
Let $p$ be a prime and $n \geq 1$. There are infinitely many
non-equivalent $MU$-products on $P(n)$. All of them induce the
same ring structure in $\mathscr D_{\SS}$ if $p$ is odd. For
$p=2$, they induce either $\mu_n$ or $\mu_n^\op$. Up to
equivalence, there is a unique commutative $MU$-product for $p$
odd and no commutative $MU$-product for $p=2$.
\end{prop}

\begin{proof}
Let $P(n)$ be endowed with the product $\mu_n$ defined as above.

Consider first the case  $p=2$. Since $v_k \equiv w_k \mod
I_{k+1}$, Proposition \ref{compofc} and \cite{jw}*{Prop.\@ 2.34}
imply that $b_{P(n)}= v_n \bar v_{n-1}^\vee \otimes \bar v
_{n-1}^\vee$. From Remark \ref{57}, we know that $\mu_n^\op =
\mu_n \circ (1+ v_n Q_{n-1} \smash Q_{n-1})$, where $Q_{n-1} \in
P(n)^*_{MU_{(p)}}(P(n))$ is the Bockstein operation associated to
$\bar v_{n-1}^\vee \in D(H_n/H_n^2[1])$. By Proposition \ref{712},
there is no commutative product on $P(n)$.

For $p$ odd, any $MU$-product on $P(n)$ can be written as
\[
\mu_n \circ \prod_{i,j}(1+\alpha_{ij}Q'_i \wedge Q'_j)
\]
for dimensional reasons, where $Q'_k$ is as in the proof of
Proposition \ref{bp prod}. Similary, for $p=2$, any $MU$-product
on $P(n)$ can be written as
\[
\mu_n \circ \prod_{i,j}(1+\alpha_{ij}Q'_i \wedge Q'_j)\circ (1+
\gamma_n v_nQ_{n-1} \wedge Q_{n-1})
\]
with $\gamma_n \in \{0,1\}$. The rest of the argument is exactly as in the proof of
Proposition \ref{bp prod}.
\end{proof}

\begin{rem}
We may consider the two degenerated cases of the family $P(n)$,
$P(0)=BP$ and $P(\infty)=\hocolim P(n)=H \mathbb F_p$, the
Eilenberg--MacLane spectrum. The former was discussed above. For
the latter, our results imply easily that it carries a unique
$MU$-product, which is commutative for all $p$.
\end{rem}

\begin{prop}
Let $BP$ be endowed with a commutative $MU$-product. Then there
are infinitely many non-equivalent $MU$-products on $P(n)$ such
that the natural map $\pi_n \: BP \to P(n)$ is multiplicative.
\end{prop}
\begin{proof}
Any product on $P(n)$ is of the form $\beta \mu_n$ with $\beta \in
\Bil(H_n/H^2_n)$, where $\mu_n$ is defined as above. By Theorem
\ref{nat}, the map $\pi_n \: BP \to \beta P(n)$ is multiplicative
if and only if $P(n)_* \otimes b_{BP} = b_{BP}^{\beta
P(n)}=\pi_n^*(b_{\beta P(n)})$. Since $BP$ is commutative,
$b_{BP}$ is trivial. Furthermore, $\pi_n \: BP \to P(n)$ is
multiplicative, by definition of $\mu_n$, and hence
$b_{BP}^{P(n)}=0$. We then deduce from Lemma \ref{technicii} that
$\pi_n \: BP \to \beta P(n)$ is multiplicative if and only if
$\pi_n^*(\beta)=0$. We easily check that there are infinitely many
such bilinear forms $\beta$ whose associated quadratic forms are
different (see Remark \ref{remclassprod}).
\end{proof}

\subsection{A non-diagonalisable product}
We aim to construct a non-diagona\-li\-zable $MU$-ring spectrum.

Let $p=2$, $I=J_2=(w_0,w_1) \subset MU_*$, as above, and $F=MU/I$.
Clearly, $F$ is a regular quotient $MU$-module, with $F_* \cong
\mathbb F_2[x_2,x_3,\ldots]$. Let $\mu$ be the smash ring product
of the products $\nu_k$ on $MU/w_k$, $ k=0,1$, from Proposition
\ref{compofc}. Let $\bar \mu = \mu \circ (1+x_2 Q_0 \smash Q_1)$,
with $Q_k$ the Bockstein operation associated to $\bar w_k^\vee
\in D(I/I^2[1])$. We claim that the product $\bar \mu$ is not
diagonalisable.

We deduce from \cite{jw}*{Prop. 2.34} and the construction of $\mu$  that the matrix of
$b_F$ with respect to the basis $\bar w_0, \bar w_1$ of $I/I^2[1]$ is
$B=\begin{pmatrix}0&0
\\0&w_2
\end{pmatrix}$.
>From Corollary \ref{newcharbil}, we deduce that the matrix of
$b_{\bar F}$ with respect to the same basis 
is given by $\bar B=\begin{pmatrix}0&x_2
\\x_2&w_2
\end{pmatrix}$.

Now assume that there exists an invertible matrix
$A=\begin{pmatrix}a&b
\\c&d
\end{pmatrix}$ with coefficients in $F_*$ such that $A^t \bar B A
=D$ is diagonal, where $A^t$ stands for the transpose of $A$. We
deduce from the equality above that
\[
(\ast) \quad (bc+ad)x_2+cdw_2=0.
\]
Since $A$ is invertible, $\det(A)$ is a unit in $F_*$. Therefore $\det(A) = ad - bc = 1$,
and hence $(\ast)$ is equivalent to
\[
(\ast\ast) \quad (1 + 2 bc)x_2+cd w_2=0.
\]
Since $|x_2|=4$ and $|w_2|=6$, there are no coefficients in $ F_*$
satisfying $(\ast\ast)$. Hence, $\bar \mu$ is not equivalent to a
diagonal product, by Proposition \ref{caradiag}.

\subsection{Morava $K$-theory $K(n)$}
The spectra  $K(n)$ can be studied as $MU$-rings, similary as we
discussed the spectra $P(n)$ above. We adopt here the point of
view of \cite{jw}*{§5} and work over the ground rings $\wh{E}(n)$
instead. We recall the definition and the notation from there. Fix
a prime number $p$. For $n > 0$, there exists a commutative
$MU_{(p)}$-algebra $\wh E(n)$ (see \cite{rog}) with
\[
\widehat{E}(n)_* \cong \lim_k
\Z_{(p)}[v_1,\ldots,v_{n-1}][v_n,v_n^{-1}]/I_n^k,
\]
where $I_n$ is the ideal generated by the regular sequence $(v_0=p, v_1, \ldots,
v_{n-1})$. The $n$-th Morava $K$-theory $K(n)$ may be defined as the regular quotient of
$\widehat{E}(n)$ by $I_n$:
\[
K(n) = \widehat{E}(n)/ I_n \cong \widehat{E}(n)/v_0
\smash_{\widehat{E}(n)} \cdots \smash_{\widehat{E}(n)}
\widehat{E}(n)/v_{n-1}.
\]
Its coefficient ring satisfies $K(n)_*\cong \F_p[v_n,v_n^{-1}].$

The following statement can be deduced from existing literature. Our methods give an
independent and quick proof. Let $Q_i\in K(n)^*_{\widehat{E}(n)}(K(n))$ be the Bockstein
operation associated to $\bar v_i^\vee\in D(I_n/I_n^2[1])$.

\begin{prop}\label{classicalmorava}
For $p$ odd, there is precisely one $\widehat{E}(n)$-product on
$K(n)$, which is commutative. For $p=2$, there are precisely two
non-equivalent $\widehat{E}(n)$-products $\mu, \bar\mu$ on $K(n)$,
both of which are non-commutative. They are related by
\[
\bar\mu = \mu^\op= \mu \circ (1+ v_n Q_{n-1} \wedge Q_{n-1}),
\]
satisfy $b_{K(n)}=b_{\bar K(n)}=v_n \bar v_{n-1}^\vee \otimes \bar
v_{n-1}^\vee$ and induce two non-equivalent $\SS$-products on
$K(n)$.
\end{prop}

\begin{proof}
The $K(n)_*$-module $I_n/I_n^2[1]$ is free with basis $\bar
v_0,\ldots ,\bar v_{n-1}$. Let first $p$ be odd. Because $|\bar
v_i\otimes\bar v_j|<|v_n|$ for all $i,j<n$, $I_n/I_n^2[1]$ admits
only the trivial bilinear form. Hence there is exactly one
$\widehat{E}(n)$-ring structure on $K(n)$ by Theorem \ref{action},
which therefore must be commutative.

Let now $p=2$. For degree reasons again, there are exactly two
bilinear forms on $I_n/I_n^2[1]$, the trivial one and $\beta=
v_n\, \bar v_{n-1}^\vee\otimes \bar v_{n-1}^\vee$. Therefore,
there are two $\widehat{E}(n)$-products $\mu$ and $\bar\mu$ on
$K(n)$, related by the formula $\bar\mu = \beta \mu=\mu \circ (1+
v_n Q_{n-1} \wedge Q_{n-1})$.

Without loss of generality, we may suppose that $\mu$ is the
diagonal product constructed in \cite{jw}*{§5.3}, whose
characteristic bilinear form $b_{K(n)}$ is  $\beta$. Hence $\mu$
is non-commutative. As a consequence, we deduce $\bar\mu=\mu^\op$,
and so $\bar\mu$ is non-commutative either. This is confirmed by
Corollary \ref{newcharbil}, which implies that $b_{\bar{K}(n)}=
b_{K(n)}- (\beta + \beta ^t)=b_{K(n)}=\beta$.

Since $\beta $ is non-alternating, we find that $\mu$ and $\mu^\op$ are not equivalent.

For the last statement, it suffices to check that the operation $Q_{n-1}$ is non-trivial
in $\mathscr D_{\SS}$, which follows from results in \cite{nassau}.
\end{proof}

\subsection{$2$-periodic Morava $K$-theory $K_n$}

We now turn to $2$-periodic Morava $K$-theory $K_n$. In this case, we have more products
than for $K(n)$, and the situation is much more interesting.

We still fix a prime number $p$ and an integer $n > 0$. There
exists a commutative $\widehat{E}(n)$-algebra $E_n$ (see
\cite{rog}), with coefficients
\[
(E_n)_*\cong \mathbb W(\F_{p^n})[[u_1, \ldots, u_{n-1}]][u^{\pm 1}],
\]
where $\mathbb W(\F_{p^n})$ is the Witt ring on $\F_{p^n}$, $|u_i|=0$ for $1\leq i \leq
n-1$ and $|u|=2$. The homomorphism induced on coefficient rings by the unit  $\eta \:
\widehat{E}(n) \to E_n $ maps $v_i$ to $u_i u^{p^i-1}$ for $1 \leq i \leq n-1$ and
$v_{n}$ to $u^{p^{n}-1}$.

Let $H_n \subset (E_n)_*$ be the ideal generated by the regular
sequence $(u_0 = p, u_1, \ldots, u_{n-1})$. The $2$-periodic
Morava $K$-theory is defined as
\[
K_n = E_n /H_n\cong E_n/u_0 \smash_{E_n} \cdots \smash_{E_n}
E_n/u_{n-1}.
\]
Its coefficient ring satisfies $(K_n)_*\cong \F_{p^n}[u,u^{-1}].$

\begin{prop}\label{2-period morava}
There are $p^n n^2$ different $E_n$-products on $K_n$, which remain different over $\SS$.
Among the $E_n$-products, none is commutative  for $p=2$; for  $p$ odd, one is
commutative if  $n=1$ and $p^n\frac{n(n-1)}{2}$ are commutative for $n>1$.
\end{prop}

\begin{proof}
The degree zero bilinear forms
\[
H_n/H_n^2[1]\otimes_{(K_n)_*} H_n/H_n^2[1] \to (K_n)_*
\]
are in bijection with the ungraded bilinear forms
\[
H_n/H_n^2\otimes_{\F_{p^n}} H_n/H_n^2 \to \F_{p^n}.
\]
It follows that there are $p^n\cdot\dim_{\F_{p^n}}(\Bil(H_n/H_n^2))= p^n n^2$ different
$E_n$-products on $K_n$.

For $p$ odd, there is a commutative $E_n$-product (see Proposition
\ref{comprod}) on $K_n$, and the set of commutative products is in
bijection with the group $\Asym(H_n/H_n^2)$, which consists of
$p^n\frac{n(n-1)}{2}$ elements for $n>1$ and one element for
$n=1$.

Let $p=2$. Using the same arguments as in the proof of Proposition 5.1 in \cite{jw}, we
construct a diagonal product $\mu$ on $K_n$ with $b_{K_n}=u\bar u_{n-1}^\vee\otimes \bar
u_{n-1}^\vee $. Hence, by Proposition  \ref{712},  $K_n$ supports no commutative
$E_n$-product.

Different $E_n$-products on $K_n$ remain different over $\SS$,
since the canonical homomorphism $(K_n)^*_{E_n}(K_n) \to
(K_n)^*_{\SS}(K_n)$ is injective. This can be deduced from
\cite{ghmr}.
\end{proof}

\begin{cor}
Up to equivalence, there are $p^n\frac{n(n+1)}{2}$ different
$E_n$-products on $K_n$. For $p$ odd, there is a unique
commutative product on $K_n$ up to equivalence.
\end{cor}
\begin{proof}
This is straightforward from Remark \ref{remclassprod} and Proposition \ref{comprod}.
\end{proof}

\begin{rem}
Our methods do not allow to determine whether non-equiva\-lent $E_n$-products on $K_n$
induce non-equivalent $\SS$-products.

A more general problem remaining open is the classification of the
set of all $\SS$-products on $K_n$, strictly as well as up to
equivalence.
\end{rem}

\begin{prop}
Any $E_n$-product on $K_n$ is diagonalizable.
\end{prop}

\begin{proof}
Apply Proposition \ref{diagprod}.
\end{proof}

\end{document}